\documentclass[reqno,12pt]{amsart}
\usepackage[utf8]{inputenc}
\usepackage[T1]{fontenc}
\usepackage{lmodern}
\usepackage{amsmath, amssymb, amsthm}
\usepackage{mathrsfs}
\usepackage{dutchcal}
\usepackage{enumitem}
\usepackage{xurl}

\numberwithin{equation}{section}

\theoremstyle{plain}
\newtheorem{theorem}{Theorem}[section]
\newtheorem{lemma}[theorem]{Lemma}
\newtheorem{corollary}[theorem]{Corollary}
\newtheorem{proposition}[theorem]{Proposition}

\theoremstyle{definition}

\newtheorem{example}[theorem]{Example}
\newtheorem{remark}[theorem]{Remark}
\newtheorem{assumption}[theorem]{Assumption}

\allowdisplaybreaks[4]

\def\dif{{\mathord{{\rm d}}}}
\def\mR{{\mathbb R}}
\def\mE{\mathbb E}
\def\mP{\mathbb P}

\usepackage{geometry}
\usepackage{color}
\usepackage[colorlinks=true]{hyperref}
\hypersetup{
linkcolor=blue,
citecolor=red,
filecolor=blue,
urlcolor=cyan
}

\usepackage[capitalize]{cleveref}

\makeatletter
\renewcommand{\subsection}{\@startsection
{subsection}
{2}
{0mm}
{0.5\baselineskip}
{0.3\baselineskip}
{\normalfont\normalsize\raggedright}}
\makeatother

\begin{document}

\title[Euler--Maruyama error for McKean--Vlasov SDEs]
{Asymptotic Error of the Euler--Maruyama Scheme for McKean--Vlasov SDEs Driven by $\alpha$-Stable Processes}

\author[Tao Wang]{Tao Wang\textsuperscript{1}}
\author[Mingkun Ye]{Mingkun Ye\textsuperscript{2,*}}
\thanks{\textsuperscript{1}\,School of Mathematics and Statistics, Jiangsu Normal University, Xuzhou, China,
  \href{mailto:taowang@jsnu.edu.cn}{taowang@jsnu.edu.cn}.}
\thanks{\textsuperscript{2}\,School of Mathematics and Statistics, Wuhan University, Wuhan, China,
  \href{mailto:mingkunye@foxmail.com}{mingkunye@foxmail.com}.}
\thanks{\textsuperscript{*}\,Corresponding author.}

\begin{abstract}
We study the continuous-time Euler--Maruyama error for McKean--Vlasov stochastic differential equations (SDEs) driven by symmetric $\alpha$-stable noise, $1<\alpha<2$. We assume globally Lipschitz coefficients, a bounded jump coefficient, and bounded Lipschitz state derivatives and law-derivative kernels. The error normalized by $(n/\log n)^{1/\alpha}$ converges jointly with the original driver in the multivariate Skorokhod $J_1$ topology, and hence stably relative to the driver's path. The limit solves a linear SDE driven by the original noise and an independent symmetric $\alpha$-stable process. Conditional centering removes both law-derivative terms from the limiting equation. For $1\leqslant p<\alpha$, the Euler schemes based on the numerical and exact marginal laws differ by $o((\log n/n)^{1/\alpha})$ in the $L^p$ norm of the path supremum. Their normalized errors have the same limiting $p$th moment of the path supremum, which is positive if $\sigma(x_0,\delta_{x_0})\partial_x\sigma(x_0,\delta_{x_0})\ne0$. For test functions in $\mathbf{C}_b^2(\mR)$, the fixed-time weak error is $o((\log n/n)^{1/\alpha})$.

\vspace{0.5em}
{\it AMS Mathematics Subject Classification}: 60H10; 60H35; 65C30.

\vspace{0.5em}
{\it Keywords}: McKean--Vlasov SDE; Euler--Maruyama scheme; asymptotic error distribution; $\alpha$-stable process.
\end{abstract}

\maketitle

\section{Introduction}\label{sec:intro}

Fix $T>0$ and $1<\alpha<2$.  
Let $L:=(L_t)_{t\geqslant 0}$ be a symmetric $\alpha$-stable L\'evy process normalized by $\mE\mathrm e^{\mathrm i\xi L_t}=\mathrm e^{-t|\xi|^\alpha}$. We use the usual augmentation $(\mathcal F_t)$ of its natural filtration. Let $\mathscr P_1(\mR)$ be the space of probability measures on $\mR$ with finite first moment, equipped with the $1$-Wasserstein distance $\mathbb W_1$, and write $\mu_Y$ for the law of $Y$. Consider the McKean--Vlasov SDE
\begin{equation}\label{eq:original_sde}
\dif X_t = b(X_t, \mu_{X_t}) \dif t + \sigma(X_{t-},\mu_{X_t}) \dif L_t, \quad X_0 = x_0 \in \mathbb{R}.
\end{equation}
The drift $b$ and jump coefficient $\sigma$ are real-valued functions
on $\mR\times\mathscr P_1(\mR)$. The jump integrand uses the predictable state $X_{t-}$ and the deterministic law $\mu_{X_t}$, continuous in $\mathbb W_1$ by Lemma \ref{lem:wellposed}(2). At each deterministic $t>0$, absence of a
fixed jump gives $\mu_{X_{t-}}=\mu_{X_t}$.

For $0\leqslant j<n$, set $t_n=jT/n$ when $jT/n<t\leqslant (j+1)T/n$,
with $t_n=0$ at $t=0$. We use $s_n$ for the same grid map evaluated at $s$.
The continuous-time Euler scheme is
\begin{equation}\label{eq:euler_scheme}
\dif X_t^n=b(X_{t_n}^n,\mu_{X_{t_n}^n})\dif t
+\sigma(X_{t_n}^n,\mu_{X_{t_n}^n})\dif L_t,\qquad X_0^n=x_0.
\end{equation}
This nonlinear Euler scheme is defined through its marginal laws $\mu_{X_{t_n}^n}$.

For $n\geqslant2$, set $u_n=(n/\log n)^{1/\alpha}$. The normalized
error is $U_t^n=u_n(X_t^n-X_t)$. This normalization appears in Jacod's Euler-error limit theorem for L\'evy-driven equations \cite[Theorem 1.2(a)]{Jacod2004AoP}.

Freezing the state and replacing the exact marginal law by the numerical
one produce distinct Euler error terms. After normalization by $u_n$, the law-feedback terms vanish, while state freezing in the jump integral yields the independent stable noise in
Theorem~\ref{thm:main}.

\begin{assumption}\label{ass:main}
\begin{enumerate}
\item The drift $b$ satisfies (a)--(c). The constant $c_b>0$ is
chosen to bound both Lipschitz estimates in (a) and (b).
  \begin{enumerate}
  \item $b$ is jointly Lipschitz in the state and the law.
  \begin{equation*}
      |b(x,\mu) - b(y,\bar\mu)| \leqslant c_b \big(|x-y| + \mathbb{W}_1(\mu,\bar\mu)\big),
  \quad \forall x,y \in \mathbb{R},\ \mu,\bar\mu \in \mathscr{P}_1(\mathbb{R}).
  \end{equation*}
  \item $\partial_x b$ exists, is bounded, and is jointly Lipschitz continuous
    in the state and measure variables.
    \begin{equation*}
    |\partial_x b(x,\mu) - \partial_x b(y,\bar\mu)|
    \leqslant c_b\bigl(|x-y|+\mathbb W_1(\mu,\bar\mu)\bigr),
    \quad \forall x,y\in\mR,\ \mu,\bar\mu\in\mathscr{P}_1(\mR).
    \end{equation*}
  \item The bounded law-derivative kernel $\mathfrak D_\mu b(x,\mu)(z)$
    is a real function of the state $x$, the law $\mu$, and the
    evaluation point $z$ of the derivative. For integrable $\xi,\eta$
    on the same atomless probability space with expectation $\widetilde{\mE}$, put $\xi_r=\eta+r(\xi-\eta)$. For every $x\in\mR$,
    \begin{equation*}
     b(x,\mu_\xi)-b(x,\mu_\eta)
    =\int_0^1\widetilde{\mE}\left[
    \mathfrak D_{\mu}b(x,\mu_{\xi_r})(\xi_r)(\xi-\eta)
    \right]\dif r.
    \end{equation*}
    For some $c_b'>0$, the kernel satisfies the following bound for all $x,y\in\mR$, $\mu,\bar\mu\in\mathscr{P}_1(\mR)$, and $z,z'\in\mR$,
    \begin{equation*}
    |\mathfrak D_{\mu}b(x,\mu)(z)-\mathfrak D_{\mu}b(y,\bar\mu)(z')|
    \leqslant c_b'\bigl(|x-y| + \mathbb{W}_1(\mu,\bar\mu) + |z-z'|\bigr).
    \end{equation*}
  \end{enumerate}
  
\item The jump coefficient $\sigma$ is bounded on $\mR\times\mathscr P_1(\mR)$. It satisfies (1)(a)--(c) with $b$ replaced by $\sigma$, using its own Lipschitz constants and law-derivative kernel $\mathfrak D_\mu\sigma$.
\end{enumerate}
\end{assumption}

The coupling identity depends only on the joint law of $(\xi,\eta)$.
For a pair on an arbitrary probability space, adjoin an atomless auxiliary
factor and apply (1)(c) on the resulting product space.

We write $\mathbf C_b^2(\mR)$ for the twice continuously differentiable
functions whose values and first two derivatives are bounded.

\begin{example}\label{ex:law-dependent-drift}
Let $b(x,\mu)=b_0(x)+F(\int\phi(z)\mu(\dif z))$, where $b_0,F,\phi\in \mathbf C_b^2(\mR)$. Then
\begin{equation*}
    \mathfrak D_\mu b(x,\mu)(z) =F'\left(\int\phi(y)\mu(\dif y)\right)\phi'(z).
\end{equation*}
Along $\xi_r=\eta+r(\xi-\eta)$, differentiation using the majorant
$\|\phi'\|_\infty|\xi-\eta|$ gives
\[
\frac{\dif}{\dif r}F\bigl(\widetilde{\mE}\phi(\xi_r)\bigr)
=F'\bigl(\widetilde{\mE}\phi(\xi_r)\bigr)
  \widetilde{\mE}[\phi'(\xi_r)(\xi-\eta)].
\]
Integrate in $r$ to get the coupling identity.
The estimate $|\int\phi\,\dif(\mu-\bar\mu)|\leqslant\|\phi'\|_\infty\mathbb W_1(\mu,\bar\mu)$ and the bounded first and second derivatives of $F$ and $\phi$
give the required Lipschitz bounds.

For an example with law dependence in both coefficients and a nonzero
limiting error, take
\begin{equation*}
b(x,\mu)=-\lambda x+\beta\int_{\mR}\sin z\,\mu(\dif z),\qquad \sigma(x,\mu)=2+\sin x+\varepsilon\int_{\mR}\cos z\,\mu(\dif z),
\end{equation*}
where $\lambda,\beta\in\mR\setminus\{0\}$ and
$0<|\varepsilon|<1/2$. Take $x_0=0$. The kernels are
\begin{equation*}
    \mathfrak D_\mu b(x,\mu)(z)=\beta\cos z,
\qquad \mathfrak D_\mu\sigma(x,\mu)(z)=-\varepsilon\sin z.
\end{equation*}
The bounded derivatives of $\sin$ and $\cos$ give the required
Lipschitz bounds, using $\mathbb W_1$ for the two law integrals.
Here $\sigma(0,\delta_0)\partial_x\sigma(0,\delta_0)=2+\varepsilon\ne0$.
Thus the limiting path moments in Corollary \ref{cor:sharp-moment} are strictly positive.
\end{example}

Lemma \ref{lem:wellposed} gives strong existence and uniqueness
for \eqref{eq:original_sde} and \eqref{eq:euler_scheme}.

\begin{theorem}
\label{thm:main}
Under Assumption \ref{ass:main}, with $X_0^n=X_0=x_0$, the normalized
errors satisfy
\begin{equation*}
(L,U^n)\Longrightarrow(L,U)
\quad\text{in }D([0,T];\mathbb R^2)\text{ with multivariate }J_1.
\end{equation*}
Consequently, $U^n$ converges stably to $U$ in $D([0,T];\mathbb R)$
with the $J_1$ topology, relative to $\sigma(L_s:0\leqslant s\leqslant T)$.
On an extension carrying a symmetric $\alpha$-stable L\'evy process $V$
independent of $L$, use the usual augmentation of their joint natural
filtration. The limit $U$ is the unique adapted c\`adl\`ag solution of
\begin{equation}\label{eq:reduced-limit}
\begin{aligned}
U_t={}&\int_0^t\partial_xb(X_s,\mu_{X_s})U_s\dif s
+\int_0^t\partial_x\sigma(X_{s-},\mu_{X_s})U_{s-}\dif L_s\\
&-\int_0^t\sigma(X_{s-},\mu_{X_s})
\partial_x\sigma(X_{s-},\mu_{X_s})\dif V_s,
\qquad U_0=0
\end{aligned}
\end{equation}
in the class satisfying $\mE\sup_{t\leqslant T}|U_t|^p<\infty$
for every $1\leqslant p<\alpha$. Its joint law with $(L,V)$ is unique.
The L\'evy measure of $V$ is
\begin{equation}\label{eq:V-Levy-measure}
\nu_V(\dif z)=\frac{T C_\alpha^2}{\alpha|z|^{1+\alpha}}\dif z,
\end{equation}
where $C_\alpha$ is defined in \eqref{eq:stable-normalization-constant}.

\end{theorem}

For a path $y$ on $[0,T]$, write $\|y\|_\infty=\sup_{t\leqslant T}|y_t|$. We have the following corollaries.
\begin{corollary}
\label{cor:sharp-moment}
Under Assumption \ref{ass:main}, for every $1\leqslant p<\alpha$,
\begin{equation}\label{eq:path-moment-asymptotic}
\lim_{n\to\infty}
\left(\frac n{\log n}\right)^{p/\alpha}
\mE\|X^n-X\|_\infty^p
=\mE\|U\|_\infty^p<\infty.
\end{equation}
If $\sigma(x_0,\delta_{x_0})\partial_x\sigma(x_0,\delta_{x_0})\ne0$, then
$\mP(\|U\|_\infty>0)=1$ and
\begin{equation*}
\bigl(\mE\|X^n-X\|_\infty^p\bigr)^{1/p}
\sim \bigl(\mE\|U\|_\infty^p\bigr)^{1/p}(\log n/n)^{1/\alpha}.
\end{equation*}
\end{corollary}

\begin{corollary}
\label{cor:feedback-bias}
Under Assumption \ref{ass:main}, set
\begin{equation}\label{eq:law-feedback-process}
u_n\bigl[f(X_{s_n}^n,\mu_{X_{s_n}^n})-f(X_{s_n}^n,\mu_{X_{s_n}})\bigr]
\mathrel{=:}G_s^{f,n},\qquad f\in\{b,\sigma\}.
\end{equation}
The corresponding drift and jump integrals are
\begin{equation*}
D_t^n=\int_0^tG_s^{b,n}\dif s,\qquad S_t^n=\int_0^tG_s^{\sigma,n}\dif L_s.
\end{equation*}
Then, for every $1\leqslant r<\infty$,
\begin{equation*}
\mE\int_0^T\bigl(|G_s^{b,n}|^r+|G_s^{\sigma,n}|^r\bigr)\dif s
\longrightarrow0.
\end{equation*}
Consequently, $\mE\|D^n\|_\infty^r\to0$ for the same $r$, and for every
$1\leqslant p<\alpha$,
\begin{equation*}
\mE\|S^n\|_\infty^p\longrightarrow0.
\end{equation*}
For each deterministic $t\in[0,T]$ and $\varphi\in\mathbf C_b^2(\mR)$,
\begin{equation}\label{eq:smooth-weak-bias}
u_n[\mE\varphi(X_t^n)-\mE\varphi(X_t)]\longrightarrow0.
\end{equation}
\end{corollary}
\begin{corollary}
\label{cor:exact-law-comparison}
Let $Y^n$ be the Euler scheme using the exact marginal law
\begin{equation*}
\dif Y_t^n=b(Y_{t_n}^n,\mu_{X_{t_n}})\dif t
+\sigma(Y_{t_n}^n,\mu_{X_{t_n}})\dif L_t,\qquad Y_0^n=x_0,
\end{equation*}
with the same driver $L$ as $X^n$. Under Assumption
\ref{ass:main}, for every $1\leqslant p<\alpha$,
\begin{equation*}
\mE\left\|u_n(X^n-Y^n)\right\|_\infty^p\longrightarrow0.
\end{equation*}
The pair $(L,u_n(Y^n-X))$ has the joint limit in Theorem \ref{thm:main},
with the corresponding stable convergence of $u_n(Y^n-X)$
relative to $\sigma(L_s:0\leqslant s\leqslant T)$. The moment asymptotic in
Corollary \ref{cor:sharp-moment} also holds for $Y^n-X$.
\end{corollary}

\subsection*{Related work and our contribution}

For a symmetric $\alpha$-stable driver, Jacod \cite[Theorem 1.2(a)]{Jacod2004AoP}
obtained the logarithmic normalization and the independent stable
noise in the limit of the grid-sampled Euler error, jointly with the
grid-sampled driver. Wang \cite[Theorem 1]{Wang2015PureJump}
also studies the grid-sampled error for pure-jump It\^o semimartingales.
For law-independent SDEs with stable noise and bounded Lipschitz
coefficients, Mikulevi\v{c}ius and Xu \cite[Proposition 3(i)]{Mikulevivcius}
give a strong upper bound at the same scale in the $L^p$ norm of the
path supremum for $1\leqslant p<\alpha$.
Jacod and Protter \cite[Theorem 3.2]{JacodandProtter1998} transfer limits
of iterated integrals to Euler-error limits when the normalized
iterated integrals satisfy their $(\star)$ condition and converge
jointly with the driver.
For general convergence results for stochastic integrals, see Kurtz
and Protter \cite{KurtzProtter1991AoP,KurtzandProtter1991}.
Brownian limits under local Lipschitz conditions and limits for
stochastic delay equations are studied by Protter et al.\ \cite{Protter2020}
and Liu et al.\ \cite{LiuWangWu2026AAP}, respectively.

For Brownian McKean--Vlasov SDEs, Liu et al.\ \cite[Theorem 4.1]{Liu2023DCDSS}
obtain a $\sqrt n$ limit for the continuous-time Euler error with
the numerical marginal law. The limit involves state and Lions
derivatives on $\mathscr P_2$ and two independent Brownian motions,
one of which is the original driver. For Brownian particle Euler
schemes, Frikha and Song \cite[Theorems 1--3]{FrikhaSong2026} derive
strong and weak bounds, including an $O(N^{-1}+h)$ weak bound for
smooth empirical-measure functionals under their regularity assumptions,
where $N$ is the particle number and $h$ the time step.

For McKean--Vlasov SDEs with additive $\alpha$-stable noise and
H\"older drift, Huang and Yang \cite{HuangYang2021} obtain
Euler--Maruyama convergence rates and propagation of chaos.
For particle Euler schemes with spectrally positive L\'evy noise and
a Brownian component, Frikha and Li \cite[Theorem 2.11]{FrikhaLi2021}
give strong $L^1$ bounds uniform over deterministic times.
Time-discretization estimates with square-integrable jumps include
the tamed schemes of Neelima et al.\ \cite{Neelima2020Levy} and the
mean-square estimates of Zhu et al.\ \cite{ZhuZhaoGan2025}.
Cavallazzi \cite{Cavallazzi2025Wellposed} studies well-posedness for
L\'evy-driven equations under finite moment assumptions that include
orders below two. Propagation of chaos for stable noise is treated
in \cite{Cavallazzi2025Wellposed,Cavallazzi2025StableChaos}.

Theorem \ref{thm:main} establishes a joint multivariate $J_1$ limit for the
continuous-time Euler error and the original stable driver when both
coefficients depend on laws in $\mathscr P_1$.
In the jump integral, freezing the exact law is negligible by
\eqref{eq:sigma-law-freezing}, while freezing the state yields the
iterated integral whose continuous-time joint limit is proved in
Section \ref{sec:noise}. Replacing the exact law by the numerical law at $X_{s_n}^n$ gives the coefficient differences $G^{b,n}$ and $G^{\sigma,n}$ in
\eqref{eq:law-feedback-process}. Section \ref{sec:proof} identifies their
limits through the law-derivative kernels and establishes joint
uniqueness in law for the resulting linear equation. Conditional
centering then makes both law-feedback terms vanish.

The joint path limit and uniform moment bounds yield the moment
convergence in Corollary \ref{cor:sharp-moment}, with a positive limit
under the stated nondegeneracy condition. The vanishing
law-feedback integrals also give the comparison
\begin{equation*}
\bigl(\mE\|X^n-Y^n\|_\infty^p\bigr)^{1/p}=o(u_n^{-1}),
\qquad 1\leqslant p<\alpha.
\end{equation*}
By Corollary \ref{cor:exact-law-comparison}, the same joint limit and
moment convergence hold for the exact-law Euler scheme.

Unless otherwise indicated, constants denoted by $C$ or $C_p$ may change between occurrences and depend on $T,\alpha,x_0$ and the coefficient bounds (and on $p$ when relevant), but not on $n$. Dependence on a localization level $K$ or $R$ is indicated by a subscript.

\section{Preliminaries}\label{sec:prelim}
We write $D([0,T];\mR^d)$ for the space of c\`adl\`ag paths,
equipped with the Skorokhod $J_1$ topology unless otherwise stated.
Convergence $x_n\to x$ in $J_1$ means that there are increasing homeomorphisms
$\lambda_n:[0,T]\to[0,T]$ such that
\begin{equation*}
\sup_{t\leqslant T}|\lambda_n(t)-t|\longrightarrow0,
\qquad
\sup_{t\leqslant T}|x_n(\lambda_n(t))-x(t)|\longrightarrow0.
\end{equation*}
In multivariate $J_1$, all coordinates share the time change $\lambda_n$.
The product of the coordinatewise $J_1$ topologies allows a separate
time change for each coordinate. See \cite{Shiryaev2003} for background.

Stable convergence relative to $\sigma(L_s:0\leqslant s\leqslant T)$ means
that, for every bounded random variable $\zeta$ measurable with respect
to this $\sigma$-field and every bounded $J_1$-continuous function
$f:D([0,T];\mR)\to\mR$,
\begin{equation*}
\mE[\zeta f(U^n)]\longrightarrow\mE[\zeta f(U)].
\end{equation*}
The expectation on the right is on the product extension, with
$\zeta$ evaluated at the original-space coordinate.

With the normalization $\mE\mathrm e^{\mathrm i\xi L_t}=\mathrm e^{-t|\xi|^\alpha}$, the L\'evy measure of $L$ is
\begin{equation*}
\nu(\dif x) = \frac{C_\alpha}{|x|^{\alpha+1}} \dif x,
\end{equation*}
where
\begin{equation}\label{eq:stable-normalization-constant}
C_\alpha
=\frac{\alpha 2^{\alpha-1}\Gamma((1+\alpha)/2)}
{\sqrt\pi\,\Gamma(1-\alpha/2)}>0.
\end{equation}
Let $\mathcal N$ be the jump measure of $L$ and
$\tilde{\mathcal N}(\dif s,\dif x)=\mathcal N(\dif s,\dif x)-\dif s\,\nu(\dif x)$
its compensated measure. Symmetry gives the L\'evy--It\^o decomposition
\begin{equation*}
L_t=\int_0^t\int_{|x|\leqslant1}x\,\tilde{\mathcal N}(\dif s,\dif x)
+\int_0^t\int_{|x|>1}x\,\mathcal N(\dif s,\dif x).
\end{equation*}

Since $\alpha>1$, $\int_{|z|>1}|z|\nu(\dif z)<\infty$ and
$\int_{|z|>1}z\,\nu(\dif z)=0$ by symmetry. Thus we may also write
$L_t=\int_0^t\int_{\mR}z\,\tilde{\mathcal N}(\dif s,\dif z)$, with
the small jumps integrated in $L^2$ and the large jumps compensated
in $L^1$. In the next two estimates we use the standard normalization in the
proofs. A change of scale for the integrator can be absorbed into $H$.

\begin{lemma}
\label{lem:stable-clock-bound}
Let $L$ be a symmetric $\alpha$-stable L\'evy process in the given
filtration. Let $\tau$ and $\tau'$ be stopping times with
$\tau\leqslant\tau'\leqslant T$, and let $H$ be predictable with $\int_0^T|H_s|^\alpha\dif s<\infty$ almost surely. For every $\varepsilon,\eta>0$,
\begin{equation}\label{eq:stable-clock-bound}
\mP\left(\sup_{\tau\leqslant t\leqslant\tau'}\left|
\int_\tau^t H_s\dif L_s\right|>\varepsilon\right)
\leqslant\mP\left(\int_\tau^{\tau'}|H_s|^\alpha\dif s>\eta\right)
+C\eta\varepsilon^{-\alpha}.
\end{equation}
Here the constant depends only on the stable law.
\end{lemma}
\begin{proof}
Let $A_t=\int_0^t\mathbf1_{\{\tau<s\leqslant\tau'\}}|H_s|^\alpha\dif s$.
To estimate the integral on $\{A_T\leqslant\eta\}$, stop at
$\tau_\eta=\inf\{t\leqslant T:A_t>\eta\}\wedge T$, with
$\inf\varnothing=\infty$. For
$K_s=H_s\mathbf1_{\{\tau<s\leqslant\tau'\}}\mathbf1_{\{s\leqslant\tau_\eta\}}$,
continuity of $A$ gives
\begin{equation*}
\int_0^T|K_s|^\alpha\dif s=A_{\tau_\eta}=A_T\wedge\eta.
\end{equation*}
For the predictable integrand $K$, the compensated integral over
$|K_sz|>\varepsilon$ is integrable because
\[
\mE\int_0^T\int_{|K_sz|>\varepsilon}|K_sz|\nu(\dif z)\dif s
=\frac{2C_\alpha}{\alpha-1}\varepsilon^{1-\alpha}\mE(A_T\wedge\eta).
\]
Its compensator vanishes by symmetry for each $K_s$. Splitting the
fully compensated integral gives
\begin{equation*}
\begin{aligned}
\int_0^t K_s\dif L_s
&=\int_0^t\int_{|K_sz|\leqslant\varepsilon}
K_sz\,\tilde{\mathcal N}(\dif s,\dif z)
+\int_0^t\int_{|K_sz|>\varepsilon}K_sz\,\mathcal N(\dif s,\dif z)\\
&\mathrel{=:} \mathcal M_t+\int_0^t\int_{|K_sz|>\varepsilon}
K_sz\,\mathcal N(\dif s,\dif z).
\end{aligned}
\end{equation*}
The number of large jumps is
\begin{equation*}
Q_T=\int_0^T\int_{|K_sz|>\varepsilon}\mathcal N(\dif s,\dif z).
\end{equation*}
By the compensation formula and the Poisson isometry,
\begin{equation*}
\mE Q_T=\frac{2C_\alpha}{\alpha}\varepsilon^{-\alpha}\mE(A_T\wedge\eta),
\qquad
\mE|\mathcal M_T|^2=\frac{2C_\alpha}{2-\alpha}
\varepsilon^{2-\alpha}\mE(A_T\wedge\eta).
\end{equation*}
On $\{Q_T=0\}$ the integral is $\mathcal M_t$, so Markov's and Doob's
inequalities give
\begin{equation*}
\begin{aligned}
\mP\left(\sup_{t\leqslant T}\left|\int_0^t K_s\dif L_s\right|>\varepsilon\right)
&\leqslant \mP(Q_T\geqslant1)+\mP\bigl(\sup_{t\leqslant T}|\mathcal M_t|>\varepsilon\bigr)\\
&\leqslant\mE Q_T+4\varepsilon^{-2}\mE|\mathcal M_T|^2
\leqslant C\varepsilon^{-\alpha}\mE(A_T\wedge\eta).
\end{aligned}
\end{equation*}
On $\{A_T\leqslant\eta\}$ the two integrals agree for $\tau\leqslant t\leqslant\tau'$, so
\begin{equation}\label{eq:stable-truncated-clock}
\begin{aligned}
\mP\left(\sup_{\tau\leqslant t\leqslant\tau'}\left|\int_\tau^tH_s\dif L_s\right|>\varepsilon\right)
&\leqslant\mP(A_T>\eta)+\mP\left(\sup_{t\leqslant T}\left|\int_0^t K_s\dif L_s\right|>\varepsilon\right)\\
&\leqslant\mP(A_T>\eta)+C\varepsilon^{-\alpha}\mE(A_T\wedge\eta).
\end{aligned}
\end{equation}
Since $A_T\wedge\eta\leqslant\eta$, this gives \eqref{eq:stable-clock-bound}.
\end{proof}

Integrating the tail bound \eqref{eq:stable-truncated-clock} gives the maximal moment inequality below. See also \cite[Theorem 3.2]{RosinskiWoyczynski1985}.

\begin{lemma}\label{lem:stable-maximal}
Let $0\leqslant S<R\leqslant T$, let $p\in(0,\alpha)$, and let $H$ be a real-valued predictable process satisfying
\begin{equation*}
\mE\left(\int_S^R|H_s|^\alpha\dif s\right)^{p/\alpha}<\infty.
\end{equation*}
Then for a symmetric $\alpha$-stable L\'evy process $L$ in this filtration,
\begin{equation*}
\mE\left[\sup_{S\leqslant t\leqslant R}\left|
\int_S^tH_s\dif L_s\right|^p\right]
\leqslant C_{p,\alpha}\mE\left(\int_S^R|H_s|^\alpha\dif s
\right)^{p/\alpha}.
\end{equation*}
The constant may also depend on the scale of $L$.
\end{lemma}

\begin{proof}
Apply \eqref{eq:stable-truncated-clock} on $[S,R]$, extending $H$ by zero off $(S,R]$, with $\varepsilon=z$ and $\eta=z^\alpha$. Integrate against $p z^{p-1}\dif z$. The tail-integral formula and Tonelli's theorem give
\begin{equation*}
\mE\sup_{S\leqslant t\leqslant R}\left|\int_S^t H_s\dif L_s\right|^p
\leqslant \mE\left(\int_S^R|H_s|^\alpha\dif s\right)^{p/\alpha}
+Cp\,\mE\int_0^\infty z^{p-\alpha-1}
\left[\left(\int_S^R|H_s|^\alpha\dif s\right)\wedge z^\alpha\right]\dif z.
\end{equation*}
For $a>0$,
\begin{equation*}
p\int_0^\infty z^{p-\alpha-1}(a\wedge z^\alpha)\dif z
=p\int_0^{a^{1/\alpha}}z^{p-1}\dif z
+pa\int_{a^{1/\alpha}}^\infty z^{p-\alpha-1}\dif z
=\frac{\alpha}{\alpha-p}a^{p/\alpha}.
\end{equation*}
The identity also holds at $a=0$, so applying it pointwise to
$\int_S^R|H_s|^\alpha\dif s$ gives the desired bound.
\end{proof}

We also use the pathwise bound
\begin{equation*}
\left(\int_S^R|H_s|^\alpha\dif s\right)^{p/\alpha}
\leqslant (R-S)^{p/\alpha}\sup_{S\leqslant s\leqslant R}|H_s|^p.
\end{equation*}

\begin{remark}\label{rem:convergence-tools}
We use Aldous' tightness criterion \cite[Theorem 1]{Aldous1978} for
c\`adl\`ag $\mathbb R^k$-valued processes with the Euclidean norm.
For $\Xi^n$ adapted to the filtration used for $\tau$, tightness of
$\sup_{t\leqslant T}|\Xi_t^n|$ and the following condition with deterministic $r$
\begin{equation*}
\lim_{\delta\downarrow0}\limsup_{n\to\infty}
\sup_{\tau\leqslant T\text{ a stopping time},\ 0\leqslant r\leqslant\delta}
\mP\left(|\Xi^n_{(\tau+r)\wedge T}-\Xi^n_\tau|>\varepsilon\right)=0
\qquad(\varepsilon>0)
\end{equation*}
imply tightness in $D([0,T];\mathbb R^k)$ with multivariate $J_1$.

For passage to the limit, let $\Xi^{n,\mathcal T}\Rightarrow\Xi^{\mathcal T}$
in multivariate $J_1$ for each deterministic partition $\mathcal T$.
For $\mathcal T=\{0=r_0<\cdots<r_m=T\}$, write $|\mathcal T|$ for
its mesh and $s_{\mathcal T}=r_k$ on $(r_k,r_{k+1}]$. The sufficient approximation conditions are
\begin{equation*}
\lim_{|\mathcal T|\to0}\limsup_{n\to\infty}
\mP(\|\Xi^n-\Xi^{n,\mathcal T}\|_\infty>\varepsilon)=0,
\qquad
\|\Xi^{\mathcal T}-\Xi\|_\infty\xrightarrow{\mP}0,
\quad |\mathcal T|\to0.
\end{equation*}
For a bounded $1$-Lipschitz test $g$ in a $J_1$ metric bounded by the
uniform distance, split each approximation error at $\varepsilon>0$.
\begin{equation*}
\begin{aligned}
|\mE g(\Xi^n)-\mE g(\Xi)|
&\leqslant |\mE g(\Xi^{n,\mathcal T})-\mE g(\Xi^{\mathcal T})|+2\varepsilon\\
&\quad+2\|g\|_\infty\bigl[
\mP(\|\Xi^n-\Xi^{n,\mathcal T}\|_\infty>\varepsilon)
+\mP(\|\Xi^{\mathcal T}-\Xi\|_\infty>\varepsilon)\bigr].
\end{aligned}
\end{equation*}
First take $n\to\infty$ at fixed $\mathcal T$, then let $|\mathcal T|\to0$
and finally $\varepsilon\downarrow0$ to obtain
$\Xi^n\Longrightarrow\Xi$ in multivariate $J_1$.
\end{remark}

\begin{lemma}
\label{lem:fixed-marginal}
Let $E,F,H$ be Polish spaces and suppose $(B_n,Y_n)\Longrightarrow(B,Y)$ in $E\times F$, with $\mu_{B_n}=\mu_B$ for every $n$. For every Borel map $f:E\to H$,
\begin{equation*}
(B_n,Y_n,f(B_n))\Longrightarrow(B,Y,f(B))
\quad\text{in the product topology.}
\end{equation*}
In particular, joint convergence with the same random
variable $B$ implies stable convergence relative to $\sigma(B)$.
\end{lemma}
\begin{proof}
For $\Phi\in\mathbf C_b(E\times F\times H)$ and $\varepsilon>0$, choose a compact $K\subset E$ on which $f$ is continuous, with $\mu_B(K^c)<\varepsilon$, by Lusin's theorem. Extend $(b,y)\mapsto\Phi(b,y,f(b))$ from $K\times F$ to a continuous $g$ of
the same sup norm bound by Tietze's theorem. Since $B_n$ and $B$ have
the same law,
\begin{equation*}
\big|\mE\Phi(B_n,Y_n,f(B_n))-\mE\Phi(B,Y,f(B))\big|
\leqslant |\mE g(B_n,Y_n)-\mE g(B,Y)|+4\|\Phi\|_\infty\varepsilon.
\end{equation*}
First use weak convergence for the test $g$, then send
$\varepsilon\downarrow0$. For $B_n=B$, take a bounded Borel
$f:E\to\mR$ and $g\in\mathbf C_b(F)$. Apply the joint convergence to
the bounded continuous test $z\,g(y)$ with $z$ clipped to
$[-\|f\|_\infty,\|f\|_\infty]$. We get $\mE[f(B)g(Y_n)]\to\mE[f(B)g(Y)]$, as required.
\end{proof}

The lemma gives product-topology convergence, which suffices for the finite coefficient vectors used below. To include a continuous $L$-measurable path in multivariate $J_1$, use the time changes for the original vector. Uniform convergence of the added path and uniform continuity of its limit give convergence under those changes.

\begin{lemma}\label{lem:moving-evaluation}
Suppose $\Xi^n\Longrightarrow \Xi$ in $D([0,T];\mR^d)$ with
multivariate $J_1$, and for every $\varepsilon>0$,
\begin{equation*}
\lim_{\delta\downarrow 0}\limsup_{n\to\infty}
\sup_{0\leqslant r\leqslant s\leqslant T,\ s-r\leqslant\delta}
\mP(|\Xi_s^n-\Xi_r^n|>\varepsilon)=0.
\end{equation*}
Then $\Xi$ has no jump at any fixed $t\in(0,T]$.
If $(B_n,\Xi^n)\Longrightarrow(B,\Xi)$ in a product of a Polish
space and this path space, then, for deterministic $v_n\to t$
in $[0,T]$,
$(B_n,\Xi_{v_n}^n)\Longrightarrow(B,\Xi_t)$.
\end{lemma}
\begin{proof}
To estimate $\Delta\Xi_t$, choose $a_j\uparrow t$ and $b_j\downarrow t$ with $a_j<t<b_j$ at deterministic continuity times of $\Xi$. The exceptional set
$\{s\in(0,T]:\mP(\Delta\Xi_s\ne0)>0\}$ is countable: for each $k$, list the finitely
many jumps of size above $1/k$ on each path and take the union of the
atoms of their time distributions over all $k$. For $t=T$, use $b_j=T$, where evaluation is always $J_1$-continuous. Since $\Xi_{b_j}-\Xi_{a_j}\to\Delta\Xi_t$
almost surely, Fatou's lemma and the Portmanteau theorem give
\begin{equation*}
\mP(|\Delta\Xi_t|>\varepsilon)
\leqslant\liminf_{j\to\infty}\mP(|\Xi_{b_j}-\Xi_{a_j}|>\varepsilon/2)
\leqslant\liminf_{j\to\infty}\liminf_{n\to\infty}
\mP(|\Xi_{b_j}^n-\Xi_{a_j}^n|>\varepsilon/2)=0 .
\end{equation*}
The last equality follows from the increment condition as $b_j-a_j\to0$.
Now evaluate jointly at the fixed time $t$. The increment condition also gives $\Xi_{v_n}^n-\Xi_t^n\xrightarrow{\mP}0$, so Slutsky's lemma
gives the assertion at $v_n$.
\end{proof}

\begin{lemma}\label{lem:partition-map}
Suppose $y^n\to y$ in $D([0,T];\mR^d)$ with multivariate
$J_1$. Fix $\mathcal T=\{0=r_0<\cdots<r_m=T\}$ and assume that
$y$ is continuous at all its partition points. For $1\leqslant\ell\leqslant r$, fix a coordinate $j_\ell$ and coefficients $a_k^{n,\ell}\to a_k^\ell$, $0\leqslant k<m$, and define
\begin{equation*}
q_t^{n,\ell}=\sum_{k=0}^{m-1}a_k^{n,\ell}
\bigl(y_{j_\ell}^n(t\wedge r_{k+1})
-y_{j_\ell}^n(t\wedge r_k)\bigr).
\end{equation*}
Define $q^\ell$ by the same sum with $y$ and $a_k^\ell$ in place of
$y^n$ and $a_k^{n,\ell}$. Then $(y^n,q^{n,1},\ldots,q^{n,r})$ converges in multivariate $J_1$ to $(y,q^1,\ldots,q^r)$. For random inputs that converge jointly in law, the same conclusion in law holds if the limiting path is almost surely continuous at the partition points.
\end{lemma}
\begin{proof}
Choose time changes $\lambda_n$ for $y^n\to y$ in $J_1$. We modify
them to fix the partition points, so the resulting $\widetilde\lambda_n$
satisfies $\widetilde\lambda_n(t\wedge r_k)=\widetilde\lambda_n(t)\wedge r_k$.
Choose $\varepsilon_n\downarrow0$ larger than $2\|\lambda_n-\mathrm{id}\|_\infty$. For large $n$, the $\varepsilon_n$-neighborhoods of the interior $r_k$ are disjoint and
lie in $(0,T)$. On each such interval replace the graph of $\lambda_n$
by the two segments joining its endpoints to $(r_k,r_k)$. This gives an
increasing homeomorphism $\widetilde\lambda_n$ tending uniformly to the
identity. The map $\lambda_n^{-1}\circ\widetilde\lambda_n$ leaves each
modified interval invariant and equals the identity elsewhere. Thus
\begin{equation*}
\|y^n\circ\widetilde\lambda_n-y\|_\infty
\leqslant\|y^n\circ\lambda_n-y\|_\infty
+\max_{0<k<m}\ \sup_{s,t\in[r_k-\varepsilon_n,r_k+\varepsilon_n]}|y_s-y_t|
\longrightarrow0.
\end{equation*}
The maximum is taken as zero for a partition without interior points.
Since $\widetilde\lambda_n$ fixes the partition, use its commutation with
$t\mapsto t\wedge r_k$ in every sum to obtain
\begin{equation*}
\|q^{n,\ell}\circ\widetilde\lambda_n-q^\ell\|_\infty
\leqslant 2\sum_{k<m}\left[
|a_k^{n,\ell}|\,\|y^n\circ\widetilde\lambda_n-y\|_\infty
+|a_k^{n,\ell}-a_k^\ell|\,\|y\|_\infty\right]\longrightarrow0.
\end{equation*}
At almost every limiting input, the path is continuous at all partition points, so the continuous mapping theorem applies to the jointly convergent paths and coefficients.
\end{proof}

The identity below changes the state at fixed law and the law at fixed
state, using the fundamental theorem of calculus and the coupling identity, respectively.
\begin{lemma}\label{lem:transport-taylor}
Let $f:\mathbb{R}\times\mathscr P_1(\mathbb{R})\to\mathbb{R}$ be continuously differentiable in its first variable. Suppose that it admits a bounded continuous kernel $\mathfrak D_{\mu}f$ satisfying the coupling identity in Assumption \ref{ass:main}(1)(c). Let $y,\bar y\in\mathbb{R}$ and let $Z,\bar Z$ be integrable real random variables on an atomless probability space with expectation $\widetilde{\mE}$. Set $\mu=\mu_Z$, $\bar\mu=\mu_{\bar Z}$ and
$Z_r=\bar Z+r(Z-\bar Z)$. Then
\begin{equation*}
f(y,\mu)-f(\bar y,\bar\mu)
=\int_0^1\partial_x f(\bar y+r(y-\bar y),\mu)(y-\bar y)\dif r +\int_0^1\widetilde{\mE}\left[
\mathfrak D_\mu f(\bar y,\mu_{Z_r})(Z_r)(Z-\bar Z)\right]\dif r.
\end{equation*}
\end{lemma}
\begin{proof}
Add and subtract $f(\bar y,\mu)$, then apply the fundamental theorem of
calculus to the state difference and the coupling identity to the law difference.
\end{proof}

\begin{remark}
A linear functional derivative on $\mathscr P_1$ is a continuous function $\frac{\delta b}{\delta\mu}(x,\mu)(z)$, with at most linear growth in $z$, uniformly for $\mu$ in compact subsets of $\mathscr P_1(\mR)$ for each fixed $x$, satisfying
\begin{equation*}
b(x,\mu)-b(x,\bar\mu) =\int_0^1\int_{\mR}
\frac{\delta b}{\delta\mu}(x,(1-r)\bar\mu+r\mu)(z) (\mu-\bar\mu)(\dif z)\dif r.
\end{equation*}
We use the integral formulation in Cavallazzi \cite[Definition 1]{Cavallazzi2022Poisson}.
Frikha and Song \cite[p.~4]{FrikhaSong2026} give this formula after
defining the derivative along mixture directions. To get the coupling identity in Assumption \ref{ass:main}(1)(c) from this mixture formula, differentiate along the laws of the interpolated random variables. If $\partial_z(\delta b/\delta\mu)$ exists, is bounded, and satisfies the Lipschitz bound in (1)(c), take
\begin{equation*}
\mathfrak D_\mu b(x,\mu)(z)=\partial_z\frac{\delta b}{\delta\mu}(x,\mu)(z).
\end{equation*}
This derivative is invariant under adding a constant in $z$ to $\delta b/\delta\mu$. For $r\in(0,1)$ and $h\ne0$ with $r+h\in[0,1]$, put $\mu_\lambda=(1-\lambda)\mu_{\xi_r}+\lambda\mu_{\xi_{r+h}}$. The mixture formula and the fundamental theorem of calculus in $z$ give
\begin{equation*}
\begin{aligned}
\frac{b(x,\mu_{\xi_{r+h}})-b(x,\mu_{\xi_r})}{h}
&=\frac1h\int_0^1\widetilde{\mE}\left[
\frac{\delta b}{\delta\mu}(x,\mu_\lambda)(\xi_{r+h})
-\frac{\delta b}{\delta\mu}(x,\mu_\lambda)(\xi_r)\right]\dif\lambda\\
&=\int_0^1\int_0^1\widetilde{\mE}\bigl[
\mathfrak D_\mu b(x,\mu_\lambda)(\xi_{r+\theta h})(\xi-\eta)
\bigr]\dif\theta\dif\lambda .
\end{aligned}
\end{equation*}
Coupling with a mixture of $\xi_r$ and $\xi_{r+h}$ gives
$\mathbb W_1(\mu_\lambda,\mu_{\xi_r})\leqslant
|h|\widetilde{\mE}|\xi-\eta|$. The integrand is bounded in absolute
value before expectation by $\|\mathfrak D_\mu b\|_\infty|\xi-\eta|$.
Dominated convergence therefore gives
\begin{equation*}
\frac{\dif}{\dif r}b(x,\mu_{\xi_r})
=\widetilde{\mE}\left[
\mathfrak D_\mu b(x,\mu_{\xi_r})(\xi_r)(\xi-\eta)\right].
\end{equation*}
The bound and continuity of the kernel give continuity of the right-hand side up to the endpoints. The map $r\mapsto b(x,\mu_{\xi_r})$ is continuous on $[0,1]$. Integration proves the coupling identity.
\end{remark}

For the marginal law of $X$, write
\begin{equation*}
\sigma_t(x):=\sigma(x,\mu_{X_t}),\qquad \sigma_t'(x):=\partial_x\sigma(x,\mu_{X_t}).
\end{equation*}
The sup norms of $\sigma$ and $\partial_x\sigma$ are taken over both arguments.

\begin{remark}\label{rem:noise-bridge}
Define the normalized iterated integral by
\begin{equation*}
M_t^n=u_n\int_0^t(L_{s-}-L_{s_n})\dif L_s.
\end{equation*}
Freezing the state and the exact law gives the stochastic error
\begin{equation*}
N_t^n=u_n\int_0^t[\sigma_{s_n}(X_{s_n})-\sigma_s(X_{s-})]\dif L_s.
\end{equation*}
Write $a_t=\sigma_t(X_t)\sigma_t'(X_t)$. The coefficient bounds and continuity of the marginal law give a
bounded c\`adl\`ag process with
$a_{s-}=\sigma_s(X_{s-})\sigma_s'(X_{s-})$. The leading part of $N^n$
is given by
\begin{equation}\label{eq:jacod-exact-noise}
\sup_{t\leqslant T}\left|
N_t^n+\int_0^t a_{s-}\dif M_s^n
\right|\xrightarrow{\mP}0.
\end{equation}
Lemma \ref{lem:noise-remainder} proves this approximation, and Corollary \ref{cor:continuous-exact-noise} gives the joint limit of $(L,M^n,N^n)$.
\end{remark}

\begin{lemma}\label{lem:wellposed}
Under Assumption \ref{ass:main}, the SDE \eqref{eq:original_sde}
and the Euler scheme \eqref{eq:euler_scheme} have unique strong solutions. Moreover,
\begin{enumerate}
\item For all $1\leqslant p<\alpha$,
\begin{equation*}
\sup_{n\geqslant1}\mE\Bigl[\sup_{0\leqslant t\leqslant T}|X_t|^p+\sup_{0\leqslant t\leqslant T}|X_t^n|^p\Bigr]<\infty.
\end{equation*}

\item For all $1\leqslant p<\alpha$, there exists $C_p>0$ such that, for every $0\leqslant s\leqslant t\leqslant T$,
\begin{equation*}
\sup_{n\geqslant1}
\mE[|X_t-X_s|^p+|X_t^n-X_s^n|^p]\leqslant C_p (t-s)^{p/\alpha}.
\end{equation*}

\item There exists $C>0$ such that, for every $n\geqslant2$ and $t\in[0,T]$,
\begin{equation*}
\mE[|X_t^n-X_{t_n}^n|^{\alpha}\wedge 1] \leqslant C\frac{\log n}{n}, \quad \mE[|X_t^n-X_{t_n}^n|^{2}\wedge 1] \leqslant C\frac{1}{n}.
\end{equation*}

\item For all $1\leqslant p<\alpha$, there exists $C_{p,T}>0$ such that, for every $n\geqslant2$,
\begin{equation*}
\mE[\sup_{0\leqslant t\leqslant T}|X_t^n - X_t|^p]
\leqslant C_{p,T} \Bigl(\frac{\log n}{n}\Bigr)^{p/\alpha}.
\end{equation*}
\end{enumerate}
\end{lemma}
\begin{proof}
\emph{Strong solutions.}
Fix $1\leqslant p<\alpha$. On $[0,h]$ let $\Gamma Z$ be the right-hand
side of the integral equation \eqref{eq:original_sde} evaluated at $Z$.
We use adapted c\`adl\`ag paths starting at $x_0$ with finite norm
$(\mE\sup_{s\leqslant h}|Z_s|^p)^{1/p}$. A Cauchy sequence in this
norm has an almost surely uniformly convergent subsequence. Its
limit is adapted and c\`adl\`ag, and Fatou's lemma gives convergence
in the norm.

For $s\downarrow t$, dominated convergence gives $\mE|Z_s-Z_t|\to0$, which bounds $\mathbb W_1(\mu_{Z_s},\mu_{Z_t})$. The law is thus a deterministic Borel function of time, so the jump integrand evaluated at $Z_{s-}$ is predictable. Linear growth and
$\mathbb W_1(\mu_{Z_s},\delta_0)=\mE|Z_s|$ bound the drift moment by
$C_p(1+\mE|Z_s|^p)$. Boundedness of $\sigma$ and Lemma \ref{lem:stable-maximal} control the jump integral. Hence $\Gamma Z$
belongs to the same space.

For $\Gamma Z-\Gamma\widetilde Z$, use H\"older's inequality for the drift and Lemma \ref{lem:stable-maximal} for the jump integral. Apply the Lipschitz and coupling bounds and Jensen's inequality.
\begin{equation*}
\begin{aligned}
\mE\sup_{t\leqslant h}|(\Gamma Z)_t-(\Gamma\widetilde Z)_t|^p
&\leqslant C_p h^{p-1}\int_0^h\mE|b(Z_s,\mu_{Z_s})-b(\widetilde Z_s,\mu_{\widetilde Z_s})|^p\dif s\\
&\quad+C_p\mE\left(\int_0^h|\sigma(Z_{s-},\mu_{Z_s})-\sigma(\widetilde Z_{s-},\mu_{\widetilde Z_s})|^\alpha\dif s\right)^{p/\alpha}\\
&\leqslant C_p(h^p+h^{p/\alpha})\mE\sup_{s\leqslant h}|Z_s-\widetilde Z_s|^p.
\end{aligned}
\end{equation*}
In the path norm the Lipschitz factor is $[C_p(h^p+h^{p/\alpha})]^{1/p}$,
with $C_p$ independent of $h$. Choose $h$ so that this factor is less
than one. The same estimate holds on each successive interval with a common $\mathcal F_r$-measurable initial value at its left endpoint $r$. Use the terminal value from the previous interval, which has finite $p$th moment, to construct the strong solution on $[0,T]$.

For pathwise uniqueness, any solution $\bar X$ must belong to the contraction space. Its drift is time-integrable, and the Lipschitz bound at fixed law gives
\begin{equation*}
\int_0^T |b(0,\mu_{\bar X_s})|\dif s
\leqslant\int_0^T |b(\bar X_s,\mu_{\bar X_s})|\dif s +c_bT\|\bar X\|_\infty<\infty\quad\text{almost surely}.
\end{equation*}
The bound $|b(\bar X_s,\mu_{\bar X_s})|\leqslant c_b|\bar X_s|+|b(0,\mu_{\bar X_s})|$
and Gronwall's inequality for the path supremum give
\begin{equation*}
\|\bar X\|_\infty\leqslant \mathrm{e}^{c_bT}\left(
|x_0|+\int_0^T|b(0,\mu_{\bar X_s})|\dif s +\sup_{t\leqslant T}\left|\int_0^t\sigma(\bar X_{s-},\mu_{\bar X_s})\dif L_s\right|\right).
\end{equation*}
The time integral on the right is a finite deterministic constant, so Lemma \ref{lem:stable-maximal} gives $\mE\|\bar X\|_\infty^p<\infty$.
The contraction therefore gives pathwise uniqueness.

For $t>0$, the scheme satisfies the recursion
\begin{equation*}
X_t^n=X_{t_n}^n+b(X_{t_n}^n,\mu_{X_{t_n}^n})(t-t_n)
+\sigma(X_{t_n}^n,\mu_{X_{t_n}^n})(L_t-L_{t_n}).
\end{equation*}
Starting from $x_0$, this recursion determines $X^n$ and its laws at the partition times. On each interval the path moment of order $p$ is finite by linear growth and Lemma \ref{lem:stable-maximal}.

(1)  Linear growth gives
\begin{equation*}
|b(x,\mu)|\leqslant |b(0,\delta_0)|
+c_b\bigl(|x|+\mathbb W_1(\mu,\delta_0)\bigr).
\end{equation*}
Since $\mathbb W_1(\mu_{X_s},\delta_0)=\mE|X_s|$, Jensen's inequality yields $\mE|b(X_s,\mu_{X_s})|^p\leqslant C_p(1+\mE|X_s|^p)$. By Lemma \ref{lem:stable-maximal}, the jump integral has $p$th path moment on $[0,t]$ at most $C_p\|\sigma\|_\infty^p t^{p/\alpha}$. Combining this with H\"older's inequality for the drift gives
\begin{equation*}
\mE\sup_{r\leqslant t}|X_r|^p
\leqslant C_{p,T}\left(1+\int_0^t\mE\sup_{v\leqslant s}|X_v|^p\dif s\right).
\end{equation*}
For $X^n$, the same drift bound uses
$|X_{s_n}^n|\leqslant\sup_{v\leqslant s}|X_v^n|$ and the same uniform bound on $\sigma$. Gronwall's inequality gives both moment estimates uniformly in $n$.

(2)  For $0\leqslant s\leqslant t\leqslant T$, taking increments in the Euler scheme
and using part (1) gives
\begin{align*}
\mE\left|\int_s^t b(X_{r_n}^n,\mu_{X_{r_n}^n})\dif r\right|^p
&\leqslant(t-s)^{p-1}\int_s^t\mE|b(X_{r_n}^n,\mu_{X_{r_n}^n})|^p\dif r
\leqslant C_p(t-s)^p,\\
\mE\left|\int_s^t\sigma(X_{r_n}^n,\mu_{X_{r_n}^n})\dif L_r\right|^p
&\leqslant C_p\mE\left(\int_s^t|\sigma(X_{r_n}^n,\mu_{X_{r_n}^n})|^\alpha\dif r\right)^{p/\alpha}
\leqslant C_p(t-s)^{p/\alpha}.
\end{align*}
For $X$, use the drift $b(X_r,\mu_{X_r})$ and the predictable jump
coefficient $\sigma(X_{r-},\mu_{X_r})$ in these two bounds. In both
equations, $|x+y|^p\leqslant2^{p-1}(|x|^p+|y|^p)$ and
$(t-s)^p\leqslant C_T(t-s)^{p/\alpha}$ give (2).

  (3) For $t>0$, use the grid recursion with left endpoint $t_n$.
  The increment length satisfies $0<t-t_n\leqslant T/n$.
  For $0<h\leqslant 1\wedge T$, \eqref{eq:stable-clock-bound} with $H\equiv1$
  gives
  \begin{equation*}
  \mP(|L_h|>y)\leqslant 1\wedge(Chy^{-\alpha}),\qquad 0<y\leqslant1.
  \end{equation*}
  For $\beta\in\{\alpha,2\}$, the tail integral gives
  \begin{equation*}
  \mE(|L_h|^\beta\wedge1)
  =\int_0^1\beta y^{\beta-1}\mP(|L_h|>y)\dif y
  \leqslant h^{\beta/\alpha}+C h\int_{h^{1/\alpha}}^1 y^{\beta-\alpha-1}\dif y.
  \end{equation*}
  Hence
  \begin{equation}\label{eq:stable-truncated-increments}
  \mE(|L_h|^\alpha\wedge1)\leqslant Ch(1+|\log h|),
  \qquad
  \mE(|L_h|^2\wedge1)\leqslant Ch.
  \end{equation}
  In the grid increment, split the drift and jump terms using
 \begin{equation*}
     |x+y|^\beta\wedge1\leqslant
  2^{\beta-1}\bigl((|x|^\beta\wedge1)+(|y|^\beta\wedge1)\bigr).
 \end{equation*}
  Since $|x|^\beta\wedge1\leqslant |x|$, boundedness of $\sigma$,
  part (1) and stationarity of the increments of $L$ give
  \begin{equation*}
  \begin{aligned}
  \mE(|X_t^n-X_{t_n}^n|^\beta\wedge1)
  &\leqslant C(t-t_n)\mE|b(X_{t_n}^n,\mu_{X_{t_n}^n})|
  +C\mE(|L_t-L_{t_n}|^\beta\wedge1)\\
  &\leqslant C\bigl(t-t_n+\mE(|L_{t-t_n}|^\beta\wedge1)\bigr).
  \end{aligned}
  \end{equation*}
  For $0<t-t_n\leqslant T/n\leqslant1$, apply
  \eqref{eq:stable-truncated-increments} and monotonicity of
  $h(1+|\log h|)$ on $(0,1]$. At $t=0$ the increment is zero.
  The bound by one covers the finitely many $n$ with $T/n>1$.

  (4) Let $S<R$ be grid points. We estimate
  $\mE\sup_{S\leqslant t\leqslant R}|X_t^n-X_t|^p$ in terms of
  $\mE|X_S^n-X_S|^p$.
  In the difference of the two SDEs, take
  $(x_s^n,x_s)=(X_s^n,X_s)$ for $f=b$ and $(X_{s-}^n,X_{s-})$ for
  $f=\sigma$. Split each coefficient difference at the numerical
  value at time $s$.
  \begin{equation*}
f(X_{s_n}^n,\mu_{X_{s_n}^n})-f(x_s,\mu_{X_s})
  =[f(X_{s_n}^n,\mu_{X_{s_n}^n})-f(x_s^n,\mu_{X_s^n})] +[f(x_s^n,\mu_{X_s^n})-f(x_s,\mu_{X_s})].
\end{equation*}
  For the drift, the Lipschitz and coupling bounds and part (2) give
  \begin{align}
  \mE\left|b(X_{s_n}^n,\mu_{X_{s_n}^n})  -b(X_s^n,\mu_{X_s^n})\right|^p
  &\leqslant C\mE|X_{s_n}^n-X_s^n|^p\leqslant C(T/n)^{p/\alpha},
  \label{eq:strong-drift-grid}\\
  |b(X_s^n,\mu_{X_s^n})-b(X_s,\mu_{X_s})|
  &\leqslant c_b\bigl(|X_s^n-X_s|+\mE|X_s^n-X_s|\bigr),
  \label{eq:strong-lipschitz-parts}
  \end{align}
  The jump-coefficient difference satisfies
  \begin{equation*}
  |\sigma(X_{s-}^n,\mu_{X_s^n})-\sigma_s(X_{s-})|
  \leqslant C\bigl(|X_{s-}^n-X_{s-}|+\mE|X_s^n-X_s|\bigr).
  \end{equation*}

  For the drift integral, use H\"older's inequality on $[S,R]$,
  \begin{equation*}
  \sup_{S\leqslant t\leqslant R}\left|\int_S^t g_s\dif s\right|^p
  \leqslant (R-S)^{p-1}\int_S^R|g_s|^p\dif s.
  \end{equation*}
  Together with \eqref{eq:strong-drift-grid}--\eqref{eq:strong-lipschitz-parts}, this gives
  \begin{equation}\label{eq:strong-drift-bound}
  \begin{aligned}
  &\mE\left[\sup_{S\leqslant t\leqslant R}\left|
  \int_S^t\bigl[b(X_{s_n}^n,\mu_{X_{s_n}^n})
  -b(X_s,\mu_{X_s})\bigr]\dif s\right|^p\right]\\
  &\quad\leqslant C(R-S)^p(T/n)^{p/\alpha}
  +C(R-S)^p\mE\sup_{S\leqslant t\leqslant R}|X_t^n-X_t|^p.
  \end{aligned}
  \end{equation}

  Boundedness and Lipschitz continuity of $\sigma$ give
  \begin{equation*}
  |\sigma(x,\mu)-\sigma(y,\bar\mu)|^\alpha
  \leqslant C\bigl[(|x-y|^\alpha\wedge1)+\mathbb W_1(\mu,\bar\mu)^\alpha\bigr].
  \end{equation*}
  For the jump integral of the grid difference, apply Lemma \ref{lem:stable-maximal}
  and Jensen's inequality, using (2)--(3) and $X_{s-}^n=X_s^n$ for almost every $s$.
  \begin{equation}\label{eq:strong-noise-grid}
  \begin{aligned}
  &\mE\left[\sup_{S\leqslant t\leqslant R}\left|
    \int_S^t\bigl[\sigma(X_{s_n}^n,\mu_{X_{s_n}^n})-\sigma(X_{s-}^n,\mu_{X_s^n})\bigr]\dif L_s
    \right|^p\right]\\
  &\quad\leqslant C\left(\mE\int_S^R|\sigma(X_{s_n}^n,\mu_{X_{s_n}^n})
    -\sigma(X_{s-}^n,\mu_{X_s^n})|^\alpha\dif s\right)^{p/\alpha}\\
  &\quad\leqslant C\left(\int_S^R\left[\mE\bigl(|X_{s_n}^n-X_s^n|^\alpha\wedge1\bigr)
    +\bigl(\mE|X_{s_n}^n-X_s^n|\bigr)^\alpha\right]\dif s\right)^{p/\alpha}
  \leqslant C\left((R-S)\frac{\log n}{n}\right)^{p/\alpha}.
  \end{aligned}
  \end{equation}
  For the other jump integral, apply the same lemma and the Lipschitz bound.
  For $s>S$ the left-limit difference is bounded by the interval
  supremum, and Jensen's inequality bounds the law term. Thus
  \begin{equation}\label{eq:strong-noise-stability}
  \begin{aligned}
  &\mE\left[\sup_{S\leqslant t\leqslant R}\left|
    \int_S^t\bigl[\sigma(X_{s-}^n,\mu_{X_s^n})-\sigma_s(X_{s-})\bigr]\dif L_s
    \right|^p\right]\\
  &\;\leqslant C\mE\left(\int_S^R
    (|X_{s-}^n-X_{s-}|+\mE|X_s^n-X_s|)^\alpha\dif s\right)^{p/\alpha}
  \leqslant C(R-S)^{p/\alpha}\mE\sup_{S\leqslant t\leqslant R}|X_t^n-X_t|^p .
  \end{aligned}
  \end{equation}

  Combining \eqref{eq:strong-drift-bound}--\eqref{eq:strong-noise-stability}, and observing that
  $(T/n)^{p/\alpha}\leqslant C(\log n/n)^{p/\alpha}$, gives
  \begin{equation*}
  \begin{aligned}
  \mE\sup_{S\leqslant t\leqslant R}|X_t^n-X_t|^p
  &\leqslant C\mE|X_S^n-X_S|^p+C\left(\frac{\log n}{n}\right)^{p/\alpha}\\
  &\quad+C\bigl((R-S)^p+(R-S)^{p/\alpha}\bigr)
  \mE\sup_{S\leqslant t\leqslant R}|X_t^n-X_t|^p.
  \end{aligned}
  \end{equation*}
  Absorb the last term on intervals of length at most $h_0$ by choosing
  $h_0$ so that its coefficient is at most $1/2$, as in
  \cite[Proposition 3]{Mikulevivcius}. When $T/n\leqslant h_0/2$, group
  $\lfloor h_0 n/T\rfloor$ grid steps at a time, with a shorter final
  interval if needed. Each full interval has length in $[h_0/2,h_0]$.
  Write $0=r_0<\cdots<r_m=T$ for this partition and set
  $E_j^n=\mE\sup_{r_j\leqslant t\leqslant r_{j+1}}|X_t^n-X_t|^p$ and
  $a_n=(\log n/n)^{p/\alpha}$. Absorbing the last term and using $X_0^n=X_0$ gives
  \begin{equation*}
  E_0^n\leqslant2Ca_n,\qquad
  E_j^n\leqslant2C(E_{j-1}^n+a_n),\quad1\leqslant j<m.
  \end{equation*}
  As $m\leqslant1+2T/h_0$ uniformly in $n$, induction and summation give
  $\mE\|X^n-X\|_\infty^p\leqslant\sum_{j<m}E_j^n\leqslant C_{p,T}a_n$.
  For the finitely many $n\geqslant2$ with $T/n>h_0/2$, enlarge
  $C_{p,T}$ using the uniform path moments in (1).
\end{proof}

On the canonical path space of $L$, write the strong solution as a
Borel path functional of $L$. Its pullback to any copy of $L$ is
adapted to that copy's augmented natural filtration. The drift integral
and the predictable partition sums for the jump integral are measurable functions of the canonical path. For the latter sums, the coefficient at the left grid point tends to $\sigma_s(X_{s-})$ for each $s>0$ by the left limits of $X$ and continuity
of $\mu_{X_s}$. The bounded integrands therefore converge in
$L^\alpha(\dif s)$ on almost every path, so Lemma \ref{lem:stable-clock-bound} gives uniform convergence in probability. Transferring the integral equation by law and applying pathwise uniqueness identifies the pullback with the strong solution.

\section{Limit of the stochastic discretization error}\label{sec:noise}

The coupling $(X_t,X_r)$ gives $\mathbb W_1(\mu_{X_t},\mu_{X_r})\leqslant\mE|X_t-X_r|$. By Assumption \ref{ass:main} and Lemma \ref{lem:wellposed}(2),
\begin{equation}\label{eq:sigma-time-modulus}
\|\sigma_t-\sigma_r\|_\infty+
\|\sigma_t'-\sigma_r'\|_\infty
\leqslant C\mathbb W_1(\mu_{X_t},\mu_{X_r})
\leqslant C|t-r|^{1/\alpha}.
\end{equation}
With $X_{s_n}$ fixed, \eqref{eq:sigma-time-modulus} and Lemma \ref{lem:stable-maximal}
give, for $1\leqslant p<\alpha$,
\begin{equation}\label{eq:sigma-law-freezing}
\mE\left\|u_n\int_0^\cdot
[\sigma_{s_n}(X_{s_n})-\sigma_s(X_{s_n})]\dif L_s
\right\|_\infty^p
\leqslant C_p\left(u_n^\alpha\int_0^T(s-s_n)\dif s\right)^{p/\alpha}
\leqslant C_p(\log n)^{-p/\alpha}.
\end{equation}
\begin{lemma}\label{lem:noise-remainder}
Under Assumption \ref{ass:main}, the approximation of $N^n$ by
\begin{equation*}
-\int_0^\cdot a_{s-}\dif M_s^n
\end{equation*}
holds uniformly in probability as in \eqref{eq:jacod-exact-noise}.
\end{lemma}
\begin{proof}
Take $n$ large enough that $h=T/n\leqslant1$, and choose $1<q<\alpha$.

\emph{Step 1. Truncated increments.}
For deterministic $0\leqslant a<s\leqslant T$ with $s-a\leqslant h$,
take the left limit in the SDE and split the two integrals at
threshold $z/2$ to obtain, for $0<z\leqslant1$,
\begin{equation*}
\begin{aligned}
\mP(|X_{s-}-X_a|>z)
&\leqslant \mP\left(\int_a^s|b(X_r,\mu_{X_r})|\dif r>z/2\right)
+\mP\left(\sup_{a\leqslant v\leqslant s}
  \left|\int_a^v\sigma_r(X_{r-})\dif L_r\right|>z/2\right)\\
&\leqslant C(s-a)(z^{-1}+z^{-\alpha})
\leqslant 2C(s-a)z^{-\alpha}.
\end{aligned}
\end{equation*}
The drift bound uses Lemma \ref{lem:wellposed}(1) and Markov's
inequality. For the jump term, use \eqref{eq:stable-clock-bound} and
$\int_a^s|\sigma_r(X_{r-})|^\alpha\dif r\leqslant(s-a)\|\sigma\|_\infty^\alpha$.
Split the tail integral at $h^{1/\alpha}$ as in
\eqref{eq:stable-truncated-increments} to obtain
\begin{equation}\label{eq:exact-truncated-increment}
\begin{aligned}
\sup_{0\leqslant a<s\leqslant T,\ s-a\leqslant h}
\mE\bigl(|X_{s-}-X_a|^\alpha\wedge1\bigr)
&\leqslant \alpha\int_0^{h^{1/\alpha}}z^{\alpha-1}\dif z
+Ch\int_{h^{1/\alpha}}^1\frac{\dif z}{z}\\
&\leqslant Ch(1+|\log h|)\leqslant C_T h\log n.
\end{aligned}
\end{equation}

\emph{Step 2. Expansion between successive partition times.}
Use the predictable grid $s_n<s$ for $s>0$. Taking the left limit
in the SDE gives, after subtracting the drift and frozen jump increment,
\begin{equation*}
X_{s-}-X_{s_n}-\int_{s_n}^s b(X_r,\mu_{X_r})\dif r
-\sigma_{s_n}(X_{s_n})(L_{s-}-L_{s_n})
=\int_{s_n}^{s-}[\sigma_r(X_{r-})-\sigma_{s_n}(X_{s_n})]\dif L_r
\mathrel{=:}J_s^n .
\end{equation*}
By \eqref{eq:sigma-time-modulus} and \eqref{eq:exact-truncated-increment},
for $s_n<r<s$,
\begin{equation*}
\mE|\sigma_r(X_{r-})-\sigma_{s_n}(X_{s_n})|^\alpha
\leqslant C\mE(|X_{r-}-X_{s_n}|^\alpha\wedge1)+C(r-s_n)
\leqslant Ch\log n.
\end{equation*}
By Lemma \ref{lem:stable-maximal} and Jensen's inequality, for $0<s\leqslant T$,
\begin{equation*}
\mE|J_s^n|^q
\leqslant C_q\left(\int_{s_n}^s\mE|\sigma_r(X_{r-})-\sigma_{s_n}(X_{s_n})|^\alpha\dif r\right)^{q/\alpha}
\leqslant C(h^2\log n)^{q/\alpha}.
\end{equation*}
In the integrand of $N^n/u_n$, add and subtract $\sigma_s(X_{s_n})$.
Using the identity for $J_s^n$ gives the remainder after subtracting
the leading term.
\begin{equation*}
\begin{aligned}
&\sigma_{s_n}(X_{s_n})-\sigma_s(X_{s-})
+\sigma_s(X_{s-})\sigma_s'(X_{s-})(L_{s-}-L_{s_n})\\
&\quad=[\sigma_{s_n}(X_{s_n})-\sigma_s(X_{s_n})]
+[\sigma_s(X_{s_n})-\sigma_s(X_{s-})+\sigma_s'(X_{s-})(X_{s-}-X_{s_n})]\\
&\qquad+\sigma_s'(X_{s-})\Big(
  [\sigma_s(X_{s-})-\sigma_{s_n}(X_{s_n})](L_{s-}-L_{s_n})
  -\int_{s_n}^s b(X_r,\mu_{X_r})\dif r-J_s^n\Big)
\mathrel{=:}\rho_s^n .
\end{aligned}
\end{equation*}
By the definition of $M^n$,
\begin{equation*}
N_t^n+\int_0^t a_{s-}\dif M_s^n
=u_n\int_0^t\rho_s^n\dif L_s .
\end{equation*}
The first bracket in $\rho_s^n$ is bounded by $Ch^{1/\alpha}$ by
\eqref{eq:sigma-time-modulus}. For the second, Taylor's formula at
$X_{s-}$ and the Lipschitz bound on $\sigma_s'$ give $C|X_{s-}-X_{s_n}|^2$.
For $K>|x_0|$, stop at $\tau_K=\inf\{t\leqslant T:|X_t|\vee|L_t|>K\}$,
with $\inf\varnothing=\infty$. On $\{s\leqslant\tau_K\}$, the left limits
at $s$ and the values strictly before $s$ are bounded by $K$, so
\begin{equation*}
|L_{s-}-L_{s_n}|\vee |X_{s-}-X_{s_n}|\leqslant2K,\qquad
\int_{s_n}^s|b(X_r,\mu_{X_r})|\dif r\leqslant C_Kh.
\end{equation*}
The identity for $J_s^n$ gives $|J_s^n|\leqslant C_K$ and hence
$|X_{s-}-X_{s_n}|^2\leqslant C_K(h^2+|L_{s-}-L_{s_n}|^2+|J_s^n|)$.
By \eqref{eq:sigma-time-modulus} and Lipschitz continuity,
$|\sigma_s(X_{s-})-\sigma_{s_n}(X_{s_n})|
\leqslant C(|X_{s-}-X_{s_n}|+h^{1/\alpha})$. Use these bounds with
$2|ab|\leqslant a^2+b^2$ and $h\leqslant1$ to obtain, on
$\{s\leqslant\tau_K\}$,
\begin{equation*}
|\rho_s^n|
\leqslant C_K\bigl(h^{1/\alpha}+|X_{s-}-X_{s_n}|^2
+|L_{s-}-L_{s_n}|^2+|J_s^n|\bigr)
\leqslant C_K\bigl(h^{1/\alpha}+|L_{s-}-L_{s_n}|^2+|J_s^n|\bigr).
\end{equation*}

\emph{Step 3. The normalized remainder.}
For deterministic $s$, integrate the stable tail as in
\eqref{eq:stable-truncated-increments} with power $2\alpha$ to get
$\mE(|L_{s-}-L_{s_n}|^{2\alpha}\wedge(2K)^{2\alpha})\leqslant C_Kh$.
Also $|J_s^n|^\alpha\mathbf1_{\{s\leqslant\tau_K\}}
\leqslant C_K |J_s^n|^q$. Consequently,
\begin{equation*}
\mE\int_0^T u_n^\alpha|\rho_s^n|^\alpha
  \mathbf1_{\{s\leqslant\tau_K\}}\dif s
\leqslant C_K u_n^\alpha\bigl[h+(h^2\log n)^{q/\alpha}\bigr]
\leqslant C_K\bigl[(\log n)^{-1}
+n^{1-2q/\alpha}(\log n)^{q/\alpha-1}\bigr]
\to 0 .
\end{equation*}
The power $1-2q/\alpha$ is negative since $q>1>\alpha/2$.
Lemma \ref{lem:stable-maximal} and Jensen's inequality give
\begin{equation*}
\mE\sup_{t\leqslant T}\left|\int_0^t u_n\rho_s^n\mathbf1_{\{s\leqslant\tau_K\}}\dif L_s\right|^q
\leqslant C_q\left(\mE\int_0^T |u_n\rho_s^n|^\alpha\mathbf1_{\{s\leqslant\tau_K\}}\dif s\right)^{q/\alpha}
\longrightarrow0.
\end{equation*}
On $\{\tau_K>T\}$, the stopped and unstopped integrals agree.
By Lemmas \ref{lem:wellposed} and \ref{lem:stable-maximal},
\begin{equation*}
\limsup_n\mP\left(\sup_{t\leqslant T}\left|u_n\int_0^t\rho_s^n\dif L_s\right|>\varepsilon\right)
\leqslant \mP(\tau_K\leqslant T)\leqslant C_qK^{-q},\qquad \varepsilon>0.
\end{equation*}
Let $K\to\infty$ to obtain \eqref{eq:jacod-exact-noise}.
\end{proof}

\begin{lemma}\label{lem:continuous-noise}
For the driving process $L$ in \eqref{eq:original_sde},
\begin{equation*}
(L,M^n)\Longrightarrow(L,V)
\quad\text{in }D([0,T];\mR^2)
\end{equation*}
with the two-dimensional $J_1$ topology, where $V$ is a symmetric
$\alpha$-stable L\'evy process, independent of $L$, with L\'evy
measure \eqref{eq:V-Levy-measure}.
\end{lemma}
\begin{proof}
Put $H_s^n=u_n(L_{s-}-L_{s_n})$, so that $M_t^n=\int_0^t H_s^n\dif L_s$. For each $n$, $H^n$ is predictable and its paths are bounded on $[0,T]$. By symmetry and scaling of $\nu$, multiplication by $H_s^n$ gives the measure $|H_s^n|^\alpha\nu$
on nonzero jump sizes. The compensator of the nonzero jumps of
$M^n$ is therefore $\dif\mathcal A_s^n\,\nu(\dif z)$, where
\begin{equation}\label{eq:noise-clock}
\mathcal A_t^n
:=\int_0^t|H_s^n|^\alpha\dif s
=u_n^\alpha\int_0^t|L_{s-}-L_{s_n}|^\alpha\dif s.
\end{equation}

\emph{Step 1. Tail estimate on a unit interval.}
Write $A=\int_0^1|L_r|^\alpha\dif r$. Its tail determines the
asymptotics of $\mathcal A^n$. The identity $n\nu(n^{1/\alpha}\,\cdot)=\nu(\cdot)$ and the criterion for infinitely divisible laws recalled in \cite[Example 17, p. 262]{HultLindskog2005} give the tail measure $\nu$ for $L_1$ at scale $n^{1/\alpha}$.
On $D([0,1];\mR)\setminus\{0\}$, where $0$ denotes the zero path, the same example gives the path tail measure
\begin{equation*}
\mathfrak m(B)=\int_0^1\int_{\mR\setminus\{0\}}
\mathbf1_{\{z\mathbf1_{[r,1]}\in B\}}\nu(\dif z)\dif r.
\end{equation*}
To compute the tail of $A$, apply this measure to the functional
$\Phi(w)=\int_0^1|w(r)|^\alpha\dif r$.
A sequence converging in $J_1$ is uniformly bounded and converges
pointwise at the continuity points of its limit. Dominated
convergence gives $J_1$-continuity of $\Phi$, which is positively
homogeneous of order $\alpha$. Moreover,
\begin{equation*}
\Phi(w)\leqslant\|w\|_\infty^\alpha,\qquad
\Phi(z\mathbf1_{[r,1]})=(1-r)|z|^\alpha .
\end{equation*}
The first bound keeps $\{\Phi>1\}$ away from the zero path. For each
$r<1$, the equation $(1-r)|z|^\alpha=1$ has only two solutions in $z$,
so $\mathfrak m(\{\Phi=1\})=0$ by the density of $\nu$.
By the mapping theorem \cite[Theorem 8 and Remark 9(ii), p. 255]{HultLindskog2005},
$ n\,\mP(A>n)\to\mathfrak m(\{\Phi>1\})$, since $\Phi(n^{-1/\alpha}L)=n^{-1}A$. Monotonicity of the tail extends the limit to real $y$. Thus
\begin{equation}\label{eq:clock-tail}
\lim_{y\to\infty}y\,\mP(A>y)
=\int_0^1\nu\bigl(\{|z|>(1-r)^{-1/\alpha}\}\bigr)\dif r =\frac{2C_\alpha}{\alpha}\int_0^1(1-r)\dif r=\frac{C_\alpha}{\alpha}.
\end{equation}

\emph{Step 2. Uniform convergence of $\mathcal A^n$.}
Put $h=T/n$. Since $L_s=L_{s-}$ for Lebesgue-almost every $s$ on
each path, changing variables on a grid interval gives
\begin{equation*}
h^{-2}\int_{jh}^{(j+1)h}|L_{s-}-L_{jh}|^\alpha\dif s
=\int_0^1|h^{-1/\alpha}(L_{jh+hr}-L_{jh})|^\alpha\dif r
\mathrel{=:}A_{j,n}.
\end{equation*}
For each fixed $n$, the variables $A_{0,n},\ldots,A_{n-1,n}$ are independent and have the law of $A$ by the independent increments and scaling of $L$.

For $2\leqslant m\leqslant n$, truncate at $B_m=m\log m$. By the tail-integral formulas and \eqref{eq:clock-tail}, as $m\to\infty$,
\begin{align*}
\mE[A\mathbf1_{\{A\leqslant B_m\}}]
&=\int_0^{B_m}\mP(A>y)\dif y-B_m\mP(A>B_m)\sim \frac{C_\alpha}{\alpha}\log m,\\
\mE[A^2\mathbf1_{\{A\leqslant B_m\}}]
&\leqslant 2\int_0^{B_m} y\,\mP(A>y)\dif y=O(m\log m).
\end{align*}
Here $B_m\mP(A>B_m)=O(1)$ and $\log B_m\sim\log m$. The truncated sum divided by $B_m$ has mean tending to $C_\alpha/\alpha$. By independence,
\begin{equation*}
\operatorname{Var}\left(\frac{1}{B_m}\sum_{j<m}A_{j,n}\mathbf1_{\{A_{j,n}\leqslant B_m\}}\right)
=\frac{m}{B_m^2}\operatorname{Var}(A\mathbf1_{\{A\leqslant B_m\}})
\leqslant\frac{m}{B_m^2}\mE[A^2\mathbf1_{\{A\leqslant B_m\}}]
=O((\log m)^{-1}).
\end{equation*}
By Chebyshev's inequality and $\mP(\max_{j<m}A_{j,n}>B_m)\leqslant m\mP(A>B_m)=O((\log m)^{-1})$,
\begin{equation*}
\frac{1}{m\log m}\sum_{j=0}^{m-1}A_{j,n} \xrightarrow{\mP}\frac{C_\alpha}{\alpha},
\qquad m=m_n\to\infty,\quad m_n\leqslant n .
\end{equation*}

For fixed $t>0$, use this limit with either $m=\lfloor t/h\rfloor$ or $m=\lceil t/h\rceil$. In both cases $m\log m/(n\log n)\to t/T$, and \eqref{eq:noise-clock} gives
\begin{equation*}
\mathcal A_{mh}^n =T^2\frac{m\log m}{n\log n}\frac{1}{m\log m}\sum_{j<m}A_{j,n} \xrightarrow{\mP}\frac{TC_\alpha}{\alpha}t .
\end{equation*}
Since $\mathcal A_{\lfloor t/h\rfloor h}^n\leqslant\mathcal A_t^n
\leqslant\mathcal A_{\lceil t/h\rceil h}^n$,
we have $\mathcal A_t^n\xrightarrow{\mP}\kappa t$ for each $t$,
where $\kappa=TC_\alpha/\alpha$.
For any fixed partition $\mathcal T$ of $[0,T]$, monotonicity gives
\begin{equation*}
\|\mathcal A^n-\kappa\,\mathrm{id}\|_\infty
\leqslant\max_{r_k\in\mathcal T}|\mathcal A_{r_k}^n-\kappa r_k|+\kappa|\mathcal T|.
\end{equation*}
Let $n\to\infty$ at fixed $\mathcal T$ and then let $|\mathcal T|\to0$. We obtain
\begin{equation}\label{eq:clock-ucp}
\|\mathcal A^n-\kappa\,\mathrm{id}\|_\infty\xrightarrow{\mP}0.
\end{equation}

\emph{Step 3. Joint laws of the increments.}
Fix $u,v\in\mR$ and $0\leqslant s<t\leqslant T$. To compute the conditional characteristic function of the increment on $(s,t]$, write $uL+vM^n$ as the integral of $u+vH^n$ against $L$ and set
\begin{equation*}
\int_0^t|u+vH_r^n|^\alpha\dif r\mathrel{=:}\Psi_t^n(u,v).
\end{equation*}
By the mean value theorem,
\begin{equation*}
\bigl||u+va|^\alpha-|u|^\alpha-|va|^\alpha\bigr|
\leqslant C_{u,v}\bigl(|a|+|a|^{\alpha-1}\bigr).
\end{equation*}
Here one expands around $u$ when $|a|\leqslant1$ and around $va$ when $|a|>1$. Integrate this bound and use stable scaling for the increments at deterministic $r$, with $r-r_n\leqslant T/n$, to get
\begin{equation*}
\begin{aligned}
&\mE\sup_{t\leqslant T}\bigl|\Psi_t^n(u,v)-t|u|^\alpha-|v|^\alpha\mathcal A_t^n\bigr|
\leqslant C_{u,v}\sum_{q\in\{1,\alpha-1\}}u_n^q
\int_0^T\mE|L_{r-}-L_{r_n}|^q\dif r\\
&\quad\leqslant C_{u,v,T}\sum_{q\in\{1,\alpha-1\}}u_n^q(T/n)^{q/\alpha}
\leqslant C_{u,v,T}(\log n)^{-(\alpha-1)/\alpha}\longrightarrow0 .
\end{aligned}
\end{equation*}
With $c=|u|^\alpha+\kappa|v|^\alpha$, \eqref{eq:clock-ucp} now gives
\begin{equation}\label{eq:joint-characteristic-clock}
\sup_{t\leqslant T}|\Psi_t^n(u,v)-tc|\xrightarrow{\mP}0 .
\end{equation}

Define the compensated exponential
\begin{equation*}
\mathcal E_t^n=\mathrm e^{\mathrm i uL_t+\mathrm i vM_t^n+\Psi_t^n(u,v)} .
\end{equation*}
For $w\in\mR$, the L\'evy--Khintchine formula under full compensation gives
\begin{equation*}
\int_{\mR}\bigl(\mathrm e^{\mathrm i w z}-1-\mathrm i w z\bigr)\nu(\dif z)
=-|w|^\alpha.
\end{equation*}
The quadratic bound at zero and $\alpha>1$ at infinity make this
integral absolutely convergent. For $K>Tc$, let $\tau_K^n=\inf\{r\in[0,T]:\Psi_r^n(u,v)\geqslant K\}$, with $\inf\varnothing=\infty$. Since $\Psi^n(u,v)$ is continuous and nondecreasing, $\Psi_{T\wedge\tau_K^n}^n\leqslant K$. Using $|\mathrm e^{\mathrm i w z}-1|^2=2(1-\cos(wz))$ and $|\mathcal E_{r-}^n|=\mathrm e^{\Psi_r^n}$, the exponent identity gives
\begin{equation*}
\int_0^{T\wedge\tau_K^n}\!\int_{\mR}
\big|\mathcal E_{r-}^n(\mathrm e^{\mathrm i(u+vH_r^n)z}-1)\big|^2
\nu(\dif z)\dif r =2\int_0^{T\wedge\tau_K^n}\mathrm e^{2\Psi_r^n}\dif\Psi_r^n
=\mathrm e^{2\Psi_{T\wedge\tau_K^n}^n}-1\leqslant\mathrm e^{2K}-1.
\end{equation*}
Apply It\^o's formula with the predictable stopping indicator
$\mathbf1_{\{r\leqslant\tau_K^n\}}$. The compensated drift
$-|u+vH_r^n|^\alpha\dif r$ cancels $\dif\Psi_r^n$, giving
\begin{equation*}
\mathcal E_{t\wedge\tau_K^n}^n-1
=\int_0^{t\wedge\tau_K^n}\int_{\mR}\mathcal E_{r-}^n
\bigl(\mathrm e^{\mathrm i(u+vH_r^n)z}-1\bigr)\,
\tilde{\mathcal N}(\dif r,\dif z).
\end{equation*}
The square-integrability bound gives a martingale by the Poisson
isometry. Since the inverse of $\mathcal E_{s\wedge\tau_K^n}^n$ has modulus at most one and is $\mathcal F_s$-measurable,
\begin{equation*}
\mE\left[ \left.\frac{\mathcal E_{t\wedge\tau_K^n}^n}{\mathcal E_{s\wedge\tau_K^n}^n}
\right|\mathcal F_s\right]=1,\qquad 0\leqslant s<t\leqslant T .
\end{equation*}
Let $\zeta_n$ be uniformly bounded and $\mathcal F_s$-measurable. Use the conditional mean one of the stopped ratio and split over $\{\Psi_T^n<K\}$, where the stopping has no effect.
\begin{equation}\label{eq:noise-conditional-characteristic}
\begin{aligned}
&\left|\mE[\zeta_n\mathrm e^{\mathrm i u(L_t-L_s)+\mathrm i v(M_t^n-M_s^n)}]
-\mathrm e^{-(t-s)c}\mE\zeta_n\right|\\
&\quad\leqslant \|\zeta_n\|_{L^\infty}\mE\left|
\frac{\mathcal E_{t\wedge\tau_K^n}^n}{\mathcal E_{s\wedge\tau_K^n}^n}
-\mathrm e^{(t-s)c+\mathrm i u(L_t-L_s)+\mathrm i v(M_t^n-M_s^n)}
\right|\\
&\quad\leqslant C\mE\left[\mathbf1_{\{\Psi_T^n<K\}}
\left|\mathrm e^{\Psi_t^n-\Psi_s^n}-\mathrm e^{(t-s)c}\right|\right]
+2C\mathrm e^K\mP(\Psi_T^n\geqslant K)\longrightarrow0.
\end{aligned}
\end{equation}
On the complement both terms in the difference have modulus at most
$\mathrm e^K$. The limit follows from \eqref{eq:joint-characteristic-clock}
and boundedness. For $0=t_0<\dots<t_m\leqslant T$, apply \eqref{eq:noise-conditional-characteristic} to the last increment, taking $\zeta_n$ to be the product of the exponentials of the preceding increments. Repeat the argument backwards to obtain
\begin{equation*}
\lim_{n\to\infty}\mE\prod_{j=1}^m
\exp\bigl(\mathrm i u_j(L_{t_j}-L_{t_{j-1}})
+\mathrm i v_j(M_{t_j}^n-M_{t_{j-1}}^n)\bigr)
=\prod_{j=1}^m \mathrm e^{-(t_j-t_{j-1})(|u_j|^\alpha+\kappa|v_j|^\alpha)}.
\end{equation*}
The product gives the finite-dimensional laws of $(L,V)$, where $V$
is independent of $L$ and has L\'evy measure \eqref{eq:V-Levy-measure} given by
\begin{equation*}
\nu_V(\dif z)=\kappa\nu(\dif z)
=\frac{TC_\alpha}{\alpha}\frac{C_\alpha}{|z|^{1+\alpha}}\dif z
=\frac{TC_\alpha^2}{\alpha|z|^{1+\alpha}}\dif z.
\end{equation*}
Lemma \ref{lem:delayed-integral} uses \eqref{eq:noise-conditional-characteristic}
to show that the drivers remain L\'evy in the filtration containing the error limit.

\emph{Step 4. Joint tightness and identification of the limit.}
The path supremum of $L$ is finite almost surely. For the second
coordinate, Lemma \ref{lem:stable-clock-bound} gives
\begin{equation*}
\mP(\|M^n\|_\infty>R)
\leqslant\mP(\mathcal A_T^n>\eta)+C\eta R^{-\alpha},
\qquad R,\eta>0.
\end{equation*}
By \eqref{eq:clock-ucp}, the first term can be made uniformly small
by choosing $\eta$ large. Letting $R\to\infty$ at that fixed $\eta$
then proves tightness of the path suprema.

For a stopping time $\tau\leqslant T$ and a deterministic $0\leqslant r\leqslant\delta$, the increase of $\mathcal A^n$ between
$\tau$ and $(\tau+r)\wedge T$ is at most $2\|\mathcal A^n-\kappa\,\mathrm{id}\|_\infty+\kappa\delta$. The same lemma therefore gives
\begin{equation*}
\mP\bigl(|M^n_{(\tau+r)\wedge T}-M_\tau^n|>\varepsilon\bigr) \leqslant
\mP(2\|\mathcal A^n-\kappa\,\mathrm{id}\|_\infty+\kappa\delta>\eta)
+C\eta\varepsilon^{-\alpha}.
\end{equation*}
For $L$, the integrand is $1$, so the bound is $C\delta\varepsilon^{-\alpha}$.
For $M^n$, fix
$\eta>0$ and take $\delta$ with $\kappa\delta<\eta$. By \eqref{eq:clock-ucp}, the upper limit of the increment probability is at most $C\eta\varepsilon^{-\alpha}$ uniformly in $\tau,r$. Let $\delta\downarrow0$ and then $\eta\downarrow0$. The two coordinates use the same stopping times, so Aldous' criterion in Remark
\ref{rem:convergence-tools} gives tightness in multivariate $J_1$.

Any subsequential limit has the finite-dimensional laws from Step 3 on its dense set of deterministic continuity times and at $T$, where $J_1$ evaluation is continuous because the time changes fix $T$. These evaluations determine the c\`adl\`ag path law, so every subsequential limit has the law of $(L,V)$.
\end{proof}

\begin{corollary}\label{cor:continuous-exact-noise}
Under Assumption \ref{ass:main}, on an extension carrying $V$, set
\begin{equation*}
N_t=-\int_0^t a_{s-}\dif V_s .
\end{equation*}
Then
\begin{equation}\label{eq:continuous-exact-noise-limit}
(L,M^n,N^n)\Longrightarrow(L,V,N)
\quad\text{in }D([0,T];\mR^3)
\end{equation}
with multivariate $J_1$ topology.
\end{corollary}
\begin{proof}
By Lemma \ref{lem:noise-remainder}, it suffices to pass to the limit in
$\int_0^\cdot a_{s-}\dif M_s^n$. Fix a deterministic partition
$\mathcal T=\{0=r_0<\cdots<r_m=T\}$. By associativity of stochastic
integration and the representation of $M^n$ in Lemma \ref{lem:continuous-noise},
\begin{equation*}
\int_0^t(a_{s-}-a_{s_{\mathcal T}})\dif M_s^n
=\int_0^t(a_{s-}-a_{s_{\mathcal T}})H_s^n\dif L_s.
\end{equation*}
By \eqref{eq:clock-ucp}, every subsequence has a further one along which
the Stieltjes measures of $\mathcal A^n$ converge weakly to $\kappa\,\dif s$
almost surely. For each path and fixed $\mathcal T$, the bounded function $s\mapsto|a_{s-}-a_{s_{\mathcal T}}|^\alpha$ is continuous outside the jump times of $a$ and the partition points. These form an at most countable set of zero $\kappa\,\dif s$-measure, so weak convergence of
the Stieltjes measures gives
\begin{equation*}
\int_0^T|a_{s-}-a_{s_{\mathcal T}}|^\alpha\dif\mathcal A_s^n
\xrightarrow{\mP}
\kappa\int_0^T|a_{s-}-a_{s_{\mathcal T}}|^\alpha\dif s.
\end{equation*}
Since $\dif\mathcal A_s^n=|H_s^n|^\alpha\dif s$, Lemma \ref{lem:stable-clock-bound}
and the Portmanteau theorem give, at fixed $\mathcal T$ and for $\varepsilon,\eta>0$,
\begin{equation*}
\begin{aligned}
\limsup_{n\to\infty}\mP\left(\sup_{t\leqslant T}\left|\int_0^t(a_{s-}-a_{s_{\mathcal T}})\dif M_s^n\right|>\varepsilon\right)
&\leqslant\limsup_{n\to\infty}\mP\left(\int_0^T|a_{s-}-a_{s_{\mathcal T}}|^\alpha\dif\mathcal A_s^n>\eta\right)
+C\eta\varepsilon^{-\alpha}\\
&\leqslant\mP\left(\kappa\int_0^T|a_{s-}-a_{s_{\mathcal T}}|^\alpha\dif s\geqslant\eta\right)
+C\eta\varepsilon^{-\alpha}.
\end{aligned}
\end{equation*}
As $|\mathcal T|\to0$, $a_{s_{\mathcal T}}\to a_{s-}$ for $s>0$,
so the last probability vanishes by bounded convergence. Letting
$\eta\downarrow0$ gives the first condition in Remark \ref{rem:convergence-tools}.

At fixed $\mathcal T$, the vector $(a_{r_k})_{k<m}$ is a measurable
function of $L$. Since the $L$ marginal is fixed, Lemma \ref{lem:fixed-marginal}
adjoins this vector to the convergence in Lemma \ref{lem:continuous-noise}.
As $L$ and $V$ have no fixed jumps, Lemma \ref{lem:partition-map} gives
\begin{equation*}
\left(L,M^n,\int_0^\cdot a_{s_{\mathcal T}}\dif M_s^n\right)
\Longrightarrow
\left(L,V,\int_0^\cdot a_{s_{\mathcal T}}\dif V_s\right).
\end{equation*}
For $1\leqslant p<\alpha$, Lemma \ref{lem:stable-maximal} for $V$
and Jensen's inequality give
\begin{equation*}
\mE\sup_{t\leqslant T}\left|\int_0^t(a_{s_{\mathcal T}}-a_{s-})\dif V_s\right|^p
\leqslant C_p\left(\mE\int_0^T|a_{s_{\mathcal T}}-a_{s-}|^\alpha\dif s\right)^{p/\alpha}
\longrightarrow0,\qquad |\mathcal T|\to0.
\end{equation*}
Remark \ref{rem:convergence-tools} therefore yields
\begin{equation*}
\left(L,M^n,\int_0^\cdot a_{s-}\dif M_s^n\right)
\Longrightarrow \left(L,V,\int_0^\cdot a_{s-}\dif V_s\right)
\quad\text{in multivariate }J_1.
\end{equation*}
By \eqref{eq:jacod-exact-noise}, the last coordinate with its
sign reversed differs from $N^n$ by a term tending uniformly
to zero in probability. This proves \eqref{eq:continuous-exact-noise-limit}.
\end{proof}

\section{Error tightness and limit identification} \label{sec:proof}

\subsection{Moment bounds and path-space tightness}

Multiplying the strong error estimate in Lemma \ref{lem:wellposed}(4)
by $u_n^p$ gives the following bound.
\begin{proposition}\label{prop:unif-mom}
Under Assumption \ref{ass:main}, for all $1\leqslant p<\alpha$,
\begin{equation*}
\sup_{n\geqslant2}\mE[\sup_{t\leqslant T}|U_t^n|^p]<\infty.
\end{equation*}
\end{proposition}

At the exact law $\mu_{X_{s_n}}$, set
$\xi_s^{n,r}=X_{s_n}+r(X_{s_n}^n-X_{s_n})$ for $r\in[0,1]$ and define
\begin{equation*}
\int_0^1\partial_x f(\xi_s^{n,r},\mu_{X_{s_n}})\dif r
\mathrel{=:}A_s^{f,n},\qquad f\in\{b,\sigma\}.
\end{equation*}
These coefficients are bounded and predictable. Use \eqref{eq:law-feedback-process} for the law increment and the fundamental theorem of calculus at the fixed law $\mu_{X_{s_n}}$ for the state
increment to obtain
\begin{equation*}
\begin{aligned}
u_n\bigl[f(X_{s_n}^n,\mu_{X_{s_n}^n})-f(X_{s_n},\mu_{X_{s_n}})\bigr]
&=G_s^{f,n}+u_n\bigl[f(X_{s_n}^n,\mu_{X_{s_n}})-f(X_{s_n},\mu_{X_{s_n}})\bigr]\\
&=G_s^{f,n}+A_s^{f,n}U_{s_n}^n.
\end{aligned}
\end{equation*}
Subtracting \eqref{eq:original_sde} from \eqref{eq:euler_scheme}
and using this decomposition gives
\begin{equation}\label{eq:exact-error}
U_t^n=\int_0^t[A_s^{b,n}U_{s_n}^n+G_s^{b,n}]\dif s
+\int_0^t[A_s^{\sigma,n}U_{s_n}^n+G_s^{\sigma,n}]\dif L_s +N_t^n+R_t^{b,n}.
\end{equation}
Here $N^n$ is given in Remark \ref{rem:noise-bridge}, and the drift
remainder is
\begin{equation*}
u_n\int_0^t[b(X_{s_n},\mu_{X_{s_n}})-b(X_s,\mu_{X_s})]\dif s
\mathrel{=:}R_t^{b,n}.
\end{equation*}
\begin{lemma}\label{lem:remainder}
Under Assumption \ref{ass:main}, for every $1\leqslant q<\alpha$,
\begin{equation*}
\mE\left[\sup_{t\leqslant T}|R_t^{b,n}|^q\right]
\leqslant C_{q,T}(\log n)^{-q/\alpha}\longrightarrow0.
\end{equation*}
\end{lemma}
\begin{proof}
Use the Lipschitz and coupling bounds in the drift integral, then
H\"older's and Jensen's inequalities and the increment estimate in
Lemma \ref{lem:wellposed}(2) to get
\begin{align*}
\mE\|R^{b,n}\|_\infty^q
&\leqslant C u_n^q\mE\left(\int_0^T
 (|X_s-X_{s_n}|+\mE|X_s-X_{s_n}|)\dif s\right)^q\\
&\leqslant C_{q,T}u_n^q\int_0^T\mE|X_s-X_{s_n}|^q\dif s
\leqslant C_{q,T}u_n^q n^{-q/\alpha}
=C_{q,T}(\log n)^{-q/\alpha}.
\end{align*}
This proves the lemma.
\end{proof}
\begin{proposition}
\label{prop:error-tightness}
Under Assumption \ref{ass:main}, the family $\{(L,M^n,N^n,U^n):n\geqslant2\}$ is tight in $D([0,T];\mR^4)$ with the multivariate $J_1$ topology. In addition, for each $a>0$,
\begin{equation}\label{eq:uniform-time-continuity-error}
\lim_{\delta\downarrow0}\limsup_{n\to\infty}
\sup_{0\leqslant r\leqslant s\leqslant T,\ s-r\leqslant\delta}
\mP(|U_s^n-U_r^n|>a)=0 .
\end{equation}
\end{proposition}
\begin{proof}
Corollary \ref{cor:continuous-exact-noise} gives tight path suprema for $(L,M^n,N^n)$, and Proposition \ref{prop:unif-mom} gives the bound for $U^n$. Let $\tau_n\leqslant T$ be stopping times in the filtration of $L$, take deterministic $\delta_n\downarrow0$, and put
$\tau_n'=(\tau_n+\delta_n)\wedge T$. Write $\Delta_n Y=Y_{\tau_n'}-Y_{\tau_n}$
for these increments.

For the increments of $(L,M^n,N^n)$, use
\begin{equation*}
\left(\Delta_n L,\Delta_n M^n,\int_{\tau_n}^{\tau_n'}a_{s-}\dif M_s^n\right)
=\int_{\tau_n}^{\tau_n'}(1,H_s^n,a_{s-}H_s^n)\dif L_s.
\end{equation*}
Over $(\tau_n,\tau_n']$, the $\alpha$th powers of the three $L^\alpha$
norms are bounded by $\delta_n$, $\Delta_n\mathcal A^n$ and
$\|a\|_\infty^\alpha\Delta_n\mathcal A^n$, respectively. By \eqref{eq:clock-ucp},
$\Delta_n\mathcal A^n\leqslant2\|\mathcal A^n-\kappa\,\mathrm{id}\|_\infty
+\kappa\delta_n\xrightarrow{\mP}0$. Lemma \ref{lem:stable-clock-bound}
gives convergence to zero of the three increments. By Lemma \ref{lem:noise-remainder},
\[
\left|\Delta_n N^n+\int_{\tau_n}^{\tau_n'}a_{s-}\dif M_s^n\right|
\leqslant 2\left\|N^n+\int_0^\cdot a_{s-}\dif M_s^n\right\|_\infty
\xrightarrow{\mP}0 .
\]

For $\Delta_n U^n$, take increments in \eqref{eq:exact-error}. The
increments of $N^n$ tend to zero by the preceding bound, and the
absolute increments of $R^{b,n}$ are bounded by
$2\|R^{b,n}\|_\infty\xrightarrow{\mP}0$
by Lemma \ref{lem:remainder}. By the coupling identity,
$|G_s^{f,n}|\leqslant\|\mathfrak D_\mu f\|_\infty\mE|U_{s_n}^n|\leqslant C$.
Both remaining integrands are bounded by $C(1+\|U^n\|_\infty)$.
On $\{\|U^n\|_\infty\leqslant R\}$, the absolute drift increment is
at most $C_R\delta_n$, as is the integral of
$|A_s^{\sigma,n}U_{s_n}^n+G_s^{\sigma,n}|^\alpha$ over $(\tau_n,\tau_n']$.
Apply Lemma \ref{lem:stable-clock-bound} with $\eta=C_R\delta_n$ and use
Proposition \ref{prop:unif-mom} at $p=1$.
\begin{equation*}
\mP\left(\left|\int_{\tau_n}^{\tau_n'}(A_s^{\sigma,n} U_{s_n}^n+G_s^{\sigma,n})\dif L_s
\right|>\varepsilon\right)
\leqslant\mP(\|U^n\|_\infty>R)+C_R\delta_n\varepsilon^{-\alpha}
\leqslant C/R+C_R\delta_n\varepsilon^{-\alpha}.
\end{equation*}
For the drift, the probability of exceeding $\varepsilon$ is at most
$C/R$ once $C_R\delta_n\leqslant\varepsilon$. First send $n\to\infty$
and then $R\to\infty$. This proves $\Delta_n U^n\xrightarrow{\mP}0$.
Since $\tau_n$ and $\delta_n$ were arbitrary, Remark
\ref{rem:convergence-tools} gives multivariate $J_1$ tightness.
Taking $\tau_n$ deterministic gives \eqref{eq:uniform-time-continuity-error},
since failure of that estimate would give deterministic intervals
of vanishing length with increment probabilities bounded away from zero.
\end{proof}

To identify the state integrals along a convergent subsequence,
we use a filtration that contains the error limit and in which both
limiting drivers are L\'evy processes.
\begin{lemma}\label{lem:delayed-integral}
Let
\begin{equation*}
I_t^n=\int_0^t A_s^{\sigma,n}U_{s_n}^n\dif L_s,\qquad
K_t^n=\int_0^t A_s^{b,n}U_{s_n}^n\dif s.
\end{equation*}
Along any subsequence for which $(L,M^n,N^n,U^n)\Longrightarrow(L,V,N,Z)$ in $D([0,T];\mR^4)$ with multivariate $J_1$, the convergence extends
jointly to $I^n$ and $K^n$ in multivariate $J_1$. Their respective limits are
\begin{equation*}
I_t=\int_0^t\sigma_s'(X_{s-})Z_{s-}\dif L_s,\qquad
K_t=\int_0^t\partial_x b(X_s,\mu_{X_s})Z_s\dif s .
\end{equation*}
Here $X$ is the strong solution driven by $L$ on the limiting space.
The vector $(L,V,N,Z)$ has no fixed jumps. In the usual augmentation
of its natural filtration, $(L,V)$ is a two-dimensional L\'evy process.
In this filtration, $N$ has the integral representation in Corollary \ref{cor:continuous-exact-noise}.
\end{lemma}
\begin{proof}
We work along the subsequence in the statement.

\emph{Step 1. The limiting filtration and $N$.}
In the product-space representation of Corollary \ref{cor:continuous-exact-noise},
\begin{equation*}
\Delta N_s=-a_{s-}\Delta V_s .
\end{equation*}
Thus $L,V,N$ have no fixed jumps. The increment estimate
\eqref{eq:uniform-time-continuity-error} and Lemma \ref{lem:moving-evaluation}
give the same property for $Z$.

Let $\mathcal H^0$ be the raw natural filtration of $(L,V,N,Z)$ and
$\mathcal H$ its usual augmentation. For $0\leqslant s_1\leqslant\cdots
\leqslant s_k\leqslant r<t\leqslant T$ and $g\in\mathbf C_b(\mR^{4k})$, put
\begin{equation*}
\zeta_n=g\bigl((L,M^n,N^n,U^n)_{s_i}:1\leqslant i\leqslant k\bigr),
\qquad
\zeta=g\bigl((L,V,N,Z)_{s_i}:1\leqslant i\leqslant k\bigr).
\end{equation*}
The variables $\zeta_n$ are $\mathcal F_r$-measurable. Joint convergence
at the fixed times $s_i,r,t$ and \eqref{eq:noise-conditional-characteristic}
give, for $u,v\in\mR$,
\begin{equation*}
\mE[\zeta\mathrm e^{\mathrm i u(L_t-L_r)+\mathrm i v(V_t-V_r)}]
=\lim_{n\to\infty}\mE[\zeta_n\mathrm e^{\mathrm i u(L_t-L_r)+\mathrm i v(M_t^n-M_r^n)}]
=\mE\zeta\,\mathrm e^{-(t-r)(|u|^\alpha+\kappa|v|^\alpha)}.
\end{equation*}
The tests $\zeta$ form an algebra of cylinder functions generating
$\mathcal H_r^0$. The monotone-class theorem gives the identity for every bounded
$\mathcal H_r^0$-measurable test. It also holds after completion by
terminal null sets. For a bounded $\mathcal H_r$-measurable test,
choose $r_j\downarrow r$ with $r_j<t$. The test belongs to the completed
$\mathcal H_{r_j}^0$ for each $j$, so we can apply the identity at $r_j$.
Right continuity of $L,V$ and bounded convergence give
\begin{equation*}
\mE\left[\left.\mathrm e^{\mathrm i u(L_t-L_r)+\mathrm i v(V_t-V_r)}
\right|\mathcal H_r\right]
=\mathrm e^{-(t-r)(|u|^\alpha+\kappa|v|^\alpha)}.
\end{equation*}
The increment law depends only on $t-r$ and is independent of
$\mathcal H_r$, so $(L,V)$ is a L\'evy process in $\mathcal H$.

To identify $N$ in this filtration, fix a deterministic partition
$\mathcal T=\{0=r_0<\cdots<r_m=T\}$. The canonical $X$ is adapted to
$\mathcal H$. The finite sum
\begin{equation*}
\int_0^t a_{s_{\mathcal T}}\dif V_s
=\sum_{k<m}a_{r_k}(V_{t\wedge r_{k+1}}-V_{t\wedge r_k})
\end{equation*}
uses the same Borel function of $(L,V)$ on both spaces.
As $|\mathcal T|\to0$, $a_{s_{\mathcal T}}\to a_{s-}$ for $s>0$.
Lemma \ref{lem:stable-maximal}
at $p=1$ and Jensen's inequality give
\begin{equation*}
\mE\left\|\int_0^\cdot(a_{s_{\mathcal T}}-a_{s-})\dif V_s\right\|_\infty
\leqslant C\left(\mE\int_0^T|a_{s_{\mathcal T}}-a_{s-}|^\alpha\dif s\right)^{1/\alpha}
\longrightarrow0,\qquad |\mathcal T|\to0.
\end{equation*}
On the product space, the definition of $N$ therefore
gives $N+\int_0^\cdot a_{s_{\mathcal T}}\dif V_s\to0$ uniformly in
probability. The joint law of $(L,V,N)$ gives the same convergence on
the subsequential limit space. The estimate above identifies
$N=-\int_0^\cdot a_{s-}\dif V_s$ there.

\emph{Step 2. Approximation of $I^n$ and $K^n$.}
Use the same partition and set
\begin{equation*}
I_t^{n,\mathcal T}=\sum_{k<m}\sigma_{r_k}'(X_{r_k})U_{r_k}^n
(L_{t\wedge r_{k+1}}-L_{t\wedge r_k}),\qquad
K_t^{n,\mathcal T}=\sum_{k<m}\partial_x b(X_{r_k},\mu_{X_{r_k}})U_{r_k}^n
(t\wedge r_{k+1}-t\wedge r_k).
\end{equation*}
For the stochastic integral, the difference of the integrands gives
\begin{equation*}
\int_0^T|A_s^{\sigma,n}U_{s_n}^n-\sigma_{s_{\mathcal T}}'(X_{s_{\mathcal T}})
U_{s_{\mathcal T}}^n|^\alpha\dif s\mathrel{=:}\mathcal C_{n,\mathcal T}.
\end{equation*}
In both differences use $au-bv=(a-b)u+b(u-v)$. On
$\{\|U^n\|_\infty\leqslant R\}$, the derivative bounds give
\begin{equation*}
\begin{aligned}
\mathcal C_{n,\mathcal T}+\|K^n-K^{n,\mathcal T}\|_\infty
&\leqslant C_R\int_0^T\Bigl(
|A_s^{\sigma,n}-\sigma_{s_{\mathcal T}}'(X_{s_{\mathcal T}})|^\alpha
+|A_s^{b,n}-\partial_x b(X_{s_{\mathcal T}},\mu_{X_{s_{\mathcal T}}})|\Bigr)\dif s\\
&\quad+C\sum_{\beta\in\{1,\alpha\}}\int_0^T
|U_{s_n}^n-U_{s_{\mathcal T}}^n|^\beta\dif s .
\end{aligned}
\end{equation*}

Compare the coefficients with the derivatives at time $s$, where
$X_s=X_{s-}$ almost surely for deterministic $s$. The arguments in the
averaged derivatives satisfy
\begin{equation*}
|\xi_s^{n,r}-X_s|\leqslant |X_{s_n}^n-X_{s_n}|+|X_{s_n}-X_s|,
\qquad \mathbb W_1(\mu_{X_{s_n}},\mu_{X_s})\leqslant\mE|X_{s_n}-X_s| .
\end{equation*}
Use the Lipschitz bounds for the derivatives and
$|A_s^{\sigma,n}-\sigma_s'(X_{s-})|^\alpha\leqslant C|A_s^{\sigma,n}-\sigma_s'(X_{s-})|$
by boundedness. Lemma \ref{lem:wellposed}(2),(4) then gives
\begin{equation}\label{eq:delayed-coefficient-Lalpha}
\begin{aligned}
&\mE\int_0^T\bigl(|A_s^{\sigma,n}-\sigma_s'(X_{s-})|^\alpha
+|A_s^{b,n}-\partial_x b(X_s,\mu_{X_s})|\bigr)\dif s\\
&\quad\leqslant C\int_0^T\bigl(\mE|X_{s_n}^n-X_{s_n}|+\mE|X_{s_n}-X_s|\bigr)\dif s
\leqslant C\bigl((\log n/n)^{1/\alpha}+n^{-1/\alpha}\bigr).
\end{aligned}
\end{equation}
Since $s_{\mathcal T}<s$ for $s>0$ and $s_{\mathcal T}\to s$ as
$|\mathcal T|\to0$,
\begin{equation*}
\bigl(\sigma_{s_{\mathcal T}}'(X_{s_{\mathcal T}}),
\partial_x b(X_{s_{\mathcal T}},\mu_{X_{s_{\mathcal T}}})\bigr)
\longrightarrow
\bigl(\sigma_s'(X_{s-}),\partial_x b(X_{s-},\mu_{X_s})\bigr).
\end{equation*}
Using bounded convergence and $X_s=X_{s-}$ for almost every $s$ on each
path, the triangle inequality and \eqref{eq:delayed-coefficient-Lalpha}
give an iterated upper limit of zero for the expected integral of
the coefficient differences, first in $n$ and then as $|\mathcal T|\to0$.

For the time shifts, fix $a>0$ and $\beta\in\{1,\alpha\}$. At fixed
$\mathcal T$, for all sufficiently large $n$,
$|s_n-s_{\mathcal T}|\leqslant T/n+|\mathcal T|\leqslant 2|\mathcal T|$.
On $\{\|U^n\|_\infty\leqslant R\}$ the difference is at most $2R$.
Splitting the integral at threshold $a$ gives
\begin{equation*}
\begin{aligned}
\mE\left[\mathbf1_{\{\|U^n\|_\infty\leqslant R\}}
\int_0^T|U_{s_n}^n-U_{s_{\mathcal T}}^n|^\beta\dif s\right]
&\leqslant Ta^\beta+(2R)^\beta\int_0^T\mP(|U_{s_n}^n-U_{s_{\mathcal T}}^n|>a)\dif s\\
&\leqslant Ta^\beta+T(2R)^\beta
\sup_{0\leqslant v\leqslant w\leqslant T,\ w-v\leqslant 2|\mathcal T|}
\mP(|U_w^n-U_v^n|>a).
\end{aligned}
\end{equation*}
By \eqref{eq:uniform-time-continuity-error}, the iterated upper limit
on the left, first in $n$ and then as $|\mathcal T|\to0$, is at most
$Ta^\beta$ for every $a>0$, hence zero.

For $\eta,R>0$, split according to $\|U^n\|_\infty>R$. Proposition \ref{prop:unif-mom} and Markov's inequality give the inequality below. For fixed $R$, the preceding estimates show that the expectation term has iterated limsup zero, first as $n\to\infty$ and then as $|\mathcal T|\to0$.
\begin{equation}\label{eq:delayed-clock-approximation}
\begin{aligned}
&\limsup_{|\mathcal T|\downarrow0}\limsup_{n\to\infty}
\mP\bigl(\mathcal C_{n,\mathcal T}+\|K^n-K^{n,\mathcal T}\|_\infty>\eta\bigr)\\
&\quad\leqslant C/R+\eta^{-1}\limsup_{|\mathcal T|\downarrow0}\limsup_{n\to\infty}
\mE\left[\mathbf1_{\{\|U^n\|_\infty\leqslant R\}}
\bigl(\mathcal C_{n,\mathcal T}+\|K^n-K^{n,\mathcal T}\|_\infty\bigr)\right]\\
&\quad=C/R\longrightarrow0,\qquad R\to\infty.
\end{aligned}
\end{equation}
For $I^n$, Lemma \ref{lem:stable-clock-bound} gives
\begin{equation*}
\mP(\|I^n-I^{n,\mathcal T}\|_\infty>\varepsilon)
\leqslant\mP(\mathcal C_{n,\mathcal T}>\eta)+C\eta\varepsilon^{-\alpha}.
\end{equation*}
Taking the upper limits as in \eqref{eq:delayed-clock-approximation}
and then letting $\eta\downarrow0$ gives
\begin{equation}\label{eq:delayed-integral-approximation}
\limsup_{|\mathcal T|\downarrow0}\limsup_{n\to\infty}
\mP(\|I^n-I^{n,\mathcal T}\|_\infty>\varepsilon)=0,\qquad \varepsilon>0 .
\end{equation}

\emph{Step 3. Joint convergence.}
At fixed $\mathcal T$, Lemma \ref{lem:fixed-marginal} lets us adjoin
the state derivatives at the $r_k$, which are measurable functions
of $L$ with its fixed law. Since $Z$ has no fixed jumps, the values
$U_{r_k}^n$ converge jointly with the paths and derivative vector.
Append $t\mapsto t$ using the same time changes. Applying
Lemma \ref{lem:partition-map} with $L$ and time as integrators gives
\begin{equation*}
(L,M^n,N^n,U^n,I^{n,\mathcal T},K^{n,\mathcal T})
\Longrightarrow (L,V,N,Z,I^{Z,\mathcal T},K^{Z,\mathcal T})
\end{equation*}
in multivariate $J_1$, where the limiting sums use $Z_{r_k}$ in place
of $U_{r_k}^n$.

On the limiting space, $Z_{s_{\mathcal T}}\to Z_{s-}$ for $s>0$ as
$|\mathcal T|\to0$. Combine this with the coefficient limits above
and the bound $C\|Z\|_\infty$ to get convergence on each path in
$L^\alpha(\dif s)$ for the jump integrand and $L^1(\dif s)$ for the
drift integrand. In the time
integral, $X_s,Z_s$ agree with their left limits almost everywhere.
Lemma \ref{lem:stable-clock-bound} and direct integration give
\begin{equation*}
\|I^{Z,\mathcal T}-I\|_\infty+\|K^{Z,\mathcal T}-K\|_\infty
\xrightarrow{\mP}0,\qquad |\mathcal T|\to0.
\end{equation*}
Together with \eqref{eq:delayed-clock-approximation} and
\eqref{eq:delayed-integral-approximation}, this allows us to pass from
the sums to the joint limit of the integrals by Remark
\ref{rem:convergence-tools}.
\end{proof}

For adapted c\`adl\`ag $Z$ with $\mE\|Z\|_\infty<\infty$, use the
conventions $X_{0-}=X_0$ and $Z_{0-}=Z_0$. For $f\in\{b,\sigma\}$, define
\begin{equation*}
F_Z^f(s,x)=\mE\bigl[\mathfrak D_\mu f(x,\mu_{X_s})(X_{s-})Z_{s-}\bigr].
\end{equation*}
For each $x$, this deterministic function is left-continuous in $s$ by
the left limits of $X,Z$, the continuity of $\mu_{X_s}$ and domination
by $\|\mathfrak D_\mu f\|_\infty\|Z\|_\infty$. Together with the
spatial Lipschitz bound, this gives an adapted, left-continuous
process $F_Z^f(s,X_{s-})$ on $(0,T]$, hence a predictable integrand.

For a subsequential limit $Z$, the joint law of $(X,Z)$ determines the
functions $F_Z^f$. Using these same deterministic functions on the
original and limiting probability spaces, set
\begin{equation*}
D_t^Z=\int_0^tF_Z^b(s,X_{s-})\dif s,\qquad
S_t^Z=\int_0^tF_Z^\sigma(s,X_{s-})\dif L_s,
\end{equation*}
using the canonical $X$ and its driver $L$ on the respective space.

\begin{samepage}
\begin{lemma}\label{lem:meanfield-passage}
Suppose that, along a subsequence,
\begin{equation*}
(L,M^n,N^n,U^n)\Longrightarrow(L,V,N,Z)
\quad\text{in }D([0,T];\mR^4)\text{ with multivariate }J_1.
\end{equation*}
Along the same subsequence, jointly with $(L,M^n,N^n,U^n,I^n,K^n)$,
we have in multivariate $J_1$
\begin{equation}\label{eq:meanfield-passage}
(D^n,S^n)\Longrightarrow(D^Z,S^Z),
\end{equation}
where the pair on the right is taken on the limiting probability space.
\end{lemma}
\end{samepage}
\begin{proof}
\emph{Step 1. Convergence at fixed time and state.}
By Lemma \ref{lem:fixed-marginal}, adjoining the canonical $X$
gives joint convergence in the product of the vector's multivariate
$J_1$ topology and the $J_1$ topology for $X$.
At deterministic $s$, $X$ and $Z$ have no fixed jumps, so
$(X_s,U_s^n)\Rightarrow(X_s,Z_s)$. In addition,
$U_{s_n}^n-U_s^n\xrightarrow{\mP}0$ by \eqref{eq:uniform-time-continuity-error}, and Lemma \ref{lem:wellposed}(2),(4) gives
\begin{equation*}
\mE|X_{s_n}^n-X_s|+\mE|X_{s_n}-X_s|
\leqslant C\bigl((\log n/n)^{1/\alpha}+n^{-1/\alpha}\bigr).
\end{equation*}
Slutsky's lemma yields
\begin{equation}\label{eq:joint-copy-marginal}
(X_{s_n},X_{s_n}^n,U_{s_n}^n)\Longrightarrow(X_s,X_s,Z_s).
\end{equation}
For $1\leqslant q<\alpha$, Proposition \ref{prop:unif-mom} gives
$\sup_n\mE\|U^n\|_\infty^q<\infty$, which also bounds
$\sup_{n,s}\mE|U_{s_n}^n|^q$. Applying the Portmanteau theorem to the
$J_1$-continuous function $\|\cdot\|_\infty^q$ gives $\mE\|Z\|_\infty^q\leqslant\liminf_n\mE\|U^n\|_\infty^q$. Fix $1<q<\alpha$ for uniform integrability below.

At fixed $x$, the coupling identity expresses the law increment
through the variables in \eqref{eq:joint-copy-marginal}.
\begin{equation*}
u_n[f(x,\mu_{X_{s_n}^n})-f(x,\mu_{X_{s_n}})]
=\int_0^1\mE[\mathfrak D_\mu f(x,\mu_{\xi_s^{n,r}})(\xi_s^{n,r})
  U_{s_n}^n]\dif r
\mathrel{=:}F_n^f(s,x).
\end{equation*}
Uniformly for $r\in[0,1]$,
\begin{equation*}
\mathbb W_1(\mu_{\xi_s^{n,r}},\mu_{X_s})
\leqslant\mE|\xi_s^{n,r}-X_s|
\leqslant\mE|X_{s_n}-X_s|+\mE|X_{s_n}^n-X_{s_n}|.
\end{equation*}
For $k\geqslant1$, use Lipschitz continuity of the kernel on
$\{|U_{s_n}^n|\leqslant k\}$ and boundedness on its complement.
The coupling and $q$th moment bounds give, with $C$ independent of $k,n,s,x$,
\begin{equation*}
\begin{aligned}
\big|F_n^f(s,x)-\mE[\mathfrak D_\mu f(x,\mu_{X_s})(X_s)U_{s_n}^n]\big|
&\leqslant Ck\sup_{0\leqslant r\leqslant1}\mE|\xi_s^{n,r}-X_s|
+C\mE\bigl[|U_{s_n}^n|\mathbf1_{\{|U_{s_n}^n|>k\}}\bigr]\\
&\leqslant Ck\bigl(\mE|X_{s_n}-X_s|+\mE|X_{s_n}^n-X_{s_n}|\bigr)
+Ck^{1-q}.
\end{aligned}
\end{equation*}
In \eqref{eq:joint-copy-marginal}, replace $X_{s_n}$ by $X_s$
using $\mE|X_{s_n}-X_s|\to0$. At a fixed $s$, $(X_s,Z_s)=(X_{s-},Z_{s-})$
almost surely, with the convention at zero. Continuity of the kernel
and uniform integrability therefore give, for fixed $s,x$,
\begin{equation*}
\mE[\mathfrak D_\mu f(x,\mu_{X_s})(X_s)U_{s_n}^n]
\longrightarrow
\mE[\mathfrak D_\mu f(x,\mu_{X_s})(X_s)Z_s]=F_Z^f(s,x).
\end{equation*}
First fix $k$ and let $n\to\infty$ in these estimates to obtain
\[
\limsup_{n\to\infty}|F_n^f(s,x)-F_Z^f(s,x)|\leqslant C k^{1-q}.
\]
Then let $k\to\infty$. Since $q>1$, this gives $F_n^f(s,x)\to F_Z^f(s,x)$.

\emph{Step 2. Convergence of the integrals on the original space.}
The bound $|F_n^f(s,x)|\leqslant\|\mathfrak D_\mu f\|_\infty\mE|U_{s_n}^n|$
is uniform in $n,s,x$, as is the bound for $F_Z^f$ using $\mE\|Z\|_\infty$.
The same moments and the Lipschitz bound on $\mathfrak D_\mu f$ give a common spatial Lipschitz constant. For fixed $s$ and $R$, pointwise convergence on finite nets and the common Lipschitz bound give uniform convergence on $[-R,R]$.

By \eqref{eq:law-feedback-process}, $G_s^{f,n}=F_n^f(s,X_{s_n}^n)$. Use the spatial Lipschitz bound to compare $X_{s_n}^n$ and $X_s$, then split according to $|X_s|\leqslant R$ to obtain
\begin{equation*}
\mE|G_s^{f,n}-F_Z^f(s,X_s)|
\leqslant C\mE|X_{s_n}^n-X_s|+\sup_{|x|\leqslant R}|F_n^f(s,x)-F_Z^f(s,x)|
+2C\mP(|X_s|>R).
\end{equation*}
The compact supremum can be taken over a countable dense set, so it is measurable in $s$. Bounded convergence at fixed $R$ and $X_s=X_{s-}$ almost surely at deterministic $s$ give
\begin{equation*}
\limsup_{n\to\infty}\mE\int_0^T|G_s^{f,n}-F_Z^f(s,X_{s-})|\dif s
\leqslant 2CT\mP(\|X\|_\infty>R).
\end{equation*}
Letting $R\to\infty$ gives convergence to zero in $L^1(\dif s\otimes\mP)$.
The uniform bound on both integrands gives, for $1\leqslant r<\infty$,
\begin{equation}\label{eq:meanfield-Lr}
\mE\int_0^T|G_s^{f,n}-F_Z^f(s,X_{s-})|^r\dif s
\leqslant (2C)^{r-1}\mE\int_0^T|G_s^{f,n}-F_Z^f(s,X_{s-})|\dif s
\longrightarrow0.
\end{equation}

H\"older's inequality and \eqref{eq:meanfield-Lr} give
\begin{equation*}
\mE\|D^n-D^Z\|_\infty^r
\leqslant T^{r-1}\mE\int_0^T|G_s^{b,n}-F_Z^b(s,X_{s-})|^r\dif s
\longrightarrow0.
\end{equation*}
For $1\leqslant p<\alpha$, Lemma \ref{lem:stable-maximal} and Jensen's
inequality give
\begin{equation}\label{eq:meanfield-stochastic-passage}
\mE\|S^n-S^Z\|_\infty^p
\leqslant C_p\left(\mE\int_0^T
|G_s^{\sigma,n}-F_Z^\sigma(s,X_{s-})|^\alpha\dif s
\right)^{p/\alpha}\longrightarrow0.
\end{equation}

\emph{Step 3. Joint convergence with the state integrals.}
Use the partition in Lemma \ref{lem:delayed-integral} on both spaces.
Set
\begin{equation*}
\begin{pmatrix}D_t^{Z,\mathcal T}\\ S_t^{Z,\mathcal T}\end{pmatrix}
=
\begin{pmatrix}
\sum_{k<m}F_Z^b(r_k,X_{r_k})(t\wedge r_{k+1}-t\wedge r_k)\\
\sum_{k<m}F_Z^\sigma(r_k,X_{r_k})(L_{t\wedge r_{k+1}}-L_{t\wedge r_k})
\end{pmatrix}.
\end{equation*}
For $s>0$, $s_{\mathcal T}$ approaches $s$ from the left as $|\mathcal T|\to0$.
The spatial Lipschitz bound and left continuity of $F_Z^f(\cdot,x)$ give
\begin{equation*}
|F_Z^f(s_{\mathcal T},X_{s_{\mathcal T}})-F_Z^f(s,X_{s-})|
\leqslant C|X_{s_{\mathcal T}}-X_{s-}|
+|F_Z^f(s_{\mathcal T},X_{s-})-F_Z^f(s,X_{s-})|
\longrightarrow0,\quad |\mathcal T|\to0.
\end{equation*}
At fixed $\mathcal T$ on the original space, use Step 2 with $r=p=1$ and the triangle inequality for the first bound. For the second, apply Lemma \ref{lem:stable-maximal} at $p=1$ and Jensen's inequality to the jump term.
\begin{equation*}
\begin{aligned}
&\limsup_{n\to\infty}\mE
\bigl[\|D^n-D^{Z,\mathcal T}\|_\infty+\|S^n-S^{Z,\mathcal T}\|_\infty\bigr]\\
&\quad\leqslant\mE
\bigl[\|D^Z-D^{Z,\mathcal T}\|_\infty+\|S^Z-S^{Z,\mathcal T}\|_\infty\bigr]\\
&\quad\leqslant\mE\int_0^T|F_Z^b(s_{\mathcal T},X_{s_{\mathcal T}})-F_Z^b(s,X_{s-})|\dif s\\
&\qquad+C\left(\mE\int_0^T
|F_Z^\sigma(s_{\mathcal T},X_{s_{\mathcal T}})-F_Z^\sigma(s,X_{s-})|^\alpha\dif s\right)^{1/\alpha}
\longrightarrow0,\qquad |\mathcal T|\to0.
\end{aligned}
\end{equation*}
The partition bound holds on both spaces and vanishes by bounded convergence.

At fixed $\mathcal T$, use Lemma \ref{lem:fixed-marginal} to adjoin the state derivatives and the values $F_Z^f(r_k,X_{r_k})$, all measurable functions of $L$ with its fixed law. Then evaluate $U^n$ at the fixed continuity times $r_k$ of its limit. Lemma \ref{lem:partition-map}, with the time coordinate appended,
gives in multivariate $J_1$
\begin{equation*}
(L,M^n,N^n,U^n,I^{n,\mathcal T},K^{n,\mathcal T},D^{Z,\mathcal T},S^{Z,\mathcal T})
\Longrightarrow
(L,V,N,Z,I^{Z,\mathcal T},K^{Z,\mathcal T},D^{Z,\mathcal T},S^{Z,\mathcal T}),
\end{equation*}
where the state sums on the right use $Z_{r_k}$. Markov's inequality for the law terms and Lemma \ref{lem:delayed-integral} for the state terms give the first condition in Remark \ref{rem:convergence-tools}. Take the upper limit in $n$ at fixed $\mathcal T$, then let $|\mathcal T|\to0$. On the limiting space the four partition sums converge uniformly in probability as $|\mathcal T|\to0$, giving the second condition. This proves
\eqref{eq:meanfield-passage}.
\end{proof}

\subsection{Limit equation and completion of the proof}
By Proposition \ref{prop:error-tightness}, every subsequence has a further subsequence such that
\begin{equation*}
(L,M^n,N^n,U^n)\Longrightarrow(L,V,N,Z).
\end{equation*}
Lemma \ref{lem:meanfield-passage} gives the joint convergence in multivariate $J_1$
\begin{equation*}
(L,M^n,N^n,U^n,K^n,D^n,I^n,S^n)
\Longrightarrow (L,V,N,Z,K,D^Z,I,S^Z).
\end{equation*}
The error equation \eqref{eq:exact-error} and Lemma \ref{lem:remainder}
give
\begin{equation*}
U^n-K^n-D^n-I^n-S^n-N^n=R^{b,n}\xrightarrow{\mP}0
\quad\text{uniformly on }[0,T].
\end{equation*}
Subtraction is continuous under the common time change in multivariate $J_1$, so the limiting difference has the law of the zero path. Hence on the limiting probability space,
\begin{equation*}
Z-K-D^Z-I-S^Z-N=0\qquad\text{almost surely}.
\end{equation*}
Substituting the limiting integrals and the representation of $N$ from Lemma \ref{lem:delayed-integral} gives
\begin{equation}\label{eq:limit_sde}
\begin{aligned}
Z_t={}&\int_0^t\bigl[\partial_xb(X_s,\mu_{X_s})Z_s
+F_Z^b(s,X_{s-})\bigr]\dif s\\
&+\int_0^t\bigl[\partial_x\sigma(X_{s-},\mu_{X_s})Z_{s-}+F_Z^\sigma(s,X_{s-})\bigr]\dif L_s\\
&-\int_0^t\sigma(X_{s-},\mu_{X_s})
\partial_x\sigma(X_{s-},\mu_{X_s})\dif V_s,
\qquad Z_0=0.
\end{aligned}
\end{equation}
\begin{lemma}\label{lem:unique-limit}
Suppose Assumption \ref{ass:main} holds. On a filtered space with the usual conditions, suppose $(L,V)$ is a two-dimensional L\'evy
process relative to that filtration, with independent coordinates. Let $L$ be standard symmetric $\alpha$-stable and $V$ symmetric $\alpha$-stable with L\'evy measure \eqref{eq:V-Levy-measure}. Let $X$ be the canonical strong solution of \eqref{eq:original_sde} driven by $L$. Equation \eqref{eq:limit_sde} then has a unique adapted c\`adl\`ag solution with $\mE\|U\|_\infty^p<\infty$ for all $1\leqslant p<\alpha$. The joint law of $(L,V,U)$ is unique across such spaces.
\end{lemma}
\begin{proof}
On an interval $[r,r+h]$, let $\Gamma$ be the integral map of \eqref{eq:limit_sde} with a prescribed $\mathcal H_r$-measurable $L^p$ initial value, where $(\mathcal H_t)$ denotes the given filtration and $1\leqslant p<\alpha$. To estimate $\Gamma Z-\Gamma\widetilde Z$, keep $X,L,V$ fixed. Subtracting the expectations in the definition of $F_Z^f$ and $F_{\widetilde Z}^f$ gives, for every $x\in\mR$,
\begin{equation}\label{eq:law-kernel-lipschitz}
|F_Z^f(s,x)-F_{\widetilde Z}^f(s,x)|
\leqslant\|\mathfrak D_\mu f\|_\infty\mE|Z_{s-}-\widetilde Z_{s-}|,
\qquad f=b,\sigma.
\end{equation}
We can use left limits in the time integral since the jump times have Lebesgue measure zero. Set $\widetilde Z=0$ and $x=X_{s-}$ in \eqref{eq:law-kernel-lipschitz}. The state derivative bounds then give
\begin{equation*}
\bigl|\partial_x f(X_{s-},\mu_{X_s})Z_{s-}+F_Z^f(s,X_{s-})\bigr|
\leqslant C\bigl(|Z_{s-}|+\mE|Z_{s-}|\bigr),\qquad f=b,\sigma .
\end{equation*}
For $\beta=1,\alpha$, bound the time integral by its length times the supremum of the integrand and use Jensen's inequality.
\begin{equation*}
\mE\left(\int_r^{r+h}\bigl(|Z_{s-}|+\mE|Z_{s-}|\bigr)^\beta\dif s\right)^{p/\beta}
\leqslant C_p h^{p/\beta}\mE\sup_{r\leqslant t\leqslant r+h}|Z_t|^p.
\end{equation*}
Use $\beta=1$ for the drift and $\beta=\alpha$ with Lemma \ref{lem:stable-maximal} for the jump integral. The same lemma for $V$ bounds the forcing term,
\begin{equation*}
\mE\sup_{r\leqslant t\leqslant r+h}\left|\int_r^t
\sigma_s(X_{s-})\sigma_s'(X_{s-})\dif V_s\right|^p
\leqslant C_ph^{p/\alpha}\|\sigma\partial_x\sigma\|_\infty^p.
\end{equation*}
Thus $\Gamma$ maps the space into itself. In the difference, the initial value and the forcing cancel. Apply the two integral bounds to $Z-\widetilde Z$ by linearity to obtain
\begin{equation*}
\mE\sup_{r\leqslant t\leqslant r+h}|(\Gamma Z)_t-(\Gamma\widetilde Z)_t|^p \leqslant C_p(h^p+h^{p/\alpha})\mE\sup_{r\leqslant t\leqslant r+h}|Z_t-\widetilde Z_t|^p.
\end{equation*}
As in Lemma \ref{lem:wellposed}, the space is complete for $\bigl(\mE\sup_{r\leqslant t\leqslant r+h}|Z_t-\widetilde Z_t|^p\bigr)^{1/p}$. Choose $h$ so that the coefficient is less than one and use each terminal value to start the next interval. Uniqueness at $p=1$ identifies the constructions for all $1\leqslant p<\alpha$.

For uniqueness of the joint law, start with $U^{(0)}=0$ on the first interval and use $U^{(j+1)}=\Gamma U^{(j)}$. The contraction gives
\begin{equation*}
\mE\sup_{0\leqslant t\leqslant h}|U_t-U_t^{(j)}|^p
\leqslant [C_p(h^p+h^{p/\alpha})]^j\mE\sup_{0\leqslant t\leqslant h}|U_t|^p\longrightarrow0.
\end{equation*}
For $\mathcal T=\{0=r_0<\cdots<r_m=h\}$, the partition approximation to
$U_t^{(j+1)}$ is
\begin{align*}
&\sum_{k<m}\bigl[\partial_x b(X_{r_k},\mu_{X_{r_k}})U_{r_k}^{(j)}
+F_{U^{(j)}}^b(r_k,X_{r_k})\bigr](t\wedge r_{k+1}-t\wedge r_k)\\
&\quad+\sum_{k<m}\bigl[\sigma_{r_k}'(X_{r_k})U_{r_k}^{(j)}
+F_{U^{(j)}}^\sigma(r_k,X_{r_k})\bigr](L_{t\wedge r_{k+1}}-L_{t\wedge r_k})\\
&\quad-\sum_{k<m}\sigma_{r_k}(X_{r_k})\sigma_{r_k}'(X_{r_k})
(V_{t\wedge r_{k+1}}-V_{t\wedge r_k}).
\end{align*}
For $f=b,\sigma$, the kernel is
\begin{equation*}
F_{U^{(j)}}^f(s,x)=\mE\bigl[
\mathfrak D_\mu f(x,\mu_{X_s})(X_{s-})U_{s-}^{(j)}\bigr].
\end{equation*}
Since $X$ is a canonical function of $L$, the law of $(L,V,U^{(j)})$
determines these kernels and the joint law of the sums with both drivers. The left limits of $X,U^{(j)}$, continuity of the marginal law and left continuity of $F_{U^{(j)}}^f$ give pointwise convergence of the step integrands. Pathwise bounds give convergence in $L^1(\dif s)$ for the drift and $L^\alpha(\dif s)$ for the two jump
terms. Direct integration and Lemma \ref{lem:stable-clock-bound} make the sums converge uniformly in probability to $U^{(j+1)}$. Thus the law of $(L,V,U^{(j)})$ determines that of $(L,V,U^{(j+1)})$. Induction from $U^{(0)}=0$ and the contraction estimate determine the joint law on the first interval.

On each successive interval $[r,r+h]$, restart the iteration with
\begin{equation*}
U_t^{(0)}=U_r,\qquad U^{(j+1)}=\Gamma U^{(j)},\qquad r\leqslant t\leqslant r+h.
\end{equation*}
Apply the same partition argument jointly to
\begin{equation*}
\bigl(L,V,(U_s)_{0\leqslant s\leqslant r},U^{(j)}\bigr).
\end{equation*}
The preceding segment and its terminal value remain in this joint law,
so the Picard limit determines the law of $\bigl(L,V,(U_s)_{0\leqslant s\leqslant r+h}\bigr)$. Repeating over finitely many intervals gives the joint law on $[0,T]$.
\end{proof}

\begin{corollary}\label{cor:meanfield-vanishes}
Under Assumption \ref{ass:main}, the unique solution $U$ of \eqref{eq:limit_sde} also solves \eqref{eq:reduced-limit}. It satisfies $\mE[U_s\mid\sigma(L_r:0\leqslant r\leqslant T)]=0$ for every $s$. In particular, both measure-derivative terms in
\eqref{eq:limit_sde} vanish.
\end{corollary}
\begin{proof}
Put $\mathcal G=\sigma(L_r:0\leqslant r\leqslant T)$. Since $F_{-U}^f=-F_U^f$ for $f=b,\sigma$ and the state terms are linear in $U$, $-U$ solves \eqref{eq:limit_sde} with $(L,-V)$. This pair is a L\'evy process in the same filtration and has the law of $(L,V)$ by independence and symmetry. Lemma \ref{lem:unique-limit} therefore yields
\begin{equation*}
    (L,V,U)\overset{\mathrm d}{=}(L,-V,-U).
\end{equation*}
For every bounded Borel function $g$ of $L$, $\mE[g(L)U_s]=-\mE[g(L)U_s]=0$, so $\mE[U_s\mid\mathcal G]=0$. For $r\uparrow s$ along deterministic times, conditional expectation contracts $L^1$, so $\|\mE[U_{s-}\mid\mathcal G]\|_{L^1}\leqslant\mE|U_{s-}-U_r|\to0$ by the integrable path supremum. At zero use $U_0=0$.

For fixed $(s,x)$, the kernel evaluated at $X_{s-}$ is $\mathcal G$-measurable because $X$ is a function of $L$. Thus
\begin{equation*}
F_U^f(s,x)
=\mE\bigl[\mathfrak D_\mu f(x,\mu_{X_s})(X_{s-})\,U_{s-}\bigr]
=\mE\bigl[\mathfrak D_\mu f(x,\mu_{X_s})(X_{s-})\,\mE[U_{s-}\mid\mathcal G]\bigr]
=0,\qquad f=b,\sigma .
\end{equation*}
Substituting $F_U^f(s,X_{s-})=0$ into \eqref{eq:limit_sde} gives
\eqref{eq:reduced-limit}.
\end{proof}

\begin{proof}[Completion of the proof of Theorem \ref{thm:main}]
By Proposition \ref{prop:error-tightness}, any subsequence has a further one with a joint limit. Its error coordinate $Z$ solves \eqref{eq:limit_sde} and has $\mE\|Z\|_\infty^p<\infty$ for $1\leqslant p<\alpha$ by the moment bound in the proof of Lemma \ref{lem:meanfield-passage}. Thus Lemma \ref{lem:unique-limit} fixes its joint law with $(L,V)$, and Corollary \ref{cor:meanfield-vanishes} identifies it as a solution of \eqref{eq:reduced-limit}. Setting the law terms to zero in the contraction and Picard arguments of Lemma \ref{lem:unique-limit} gives pathwise uniqueness and uniqueness of the joint law for \eqref{eq:reduced-limit}. Thus $(L,U^n)\Longrightarrow(L,U)$ in multivariate $J_1$. The first coordinate is the original $L$, so Lemma \ref{lem:fixed-marginal} gives stable convergence relative to its path $\sigma$-field.
\end{proof}

\begin{proof}[Proof of Corollary \ref{cor:sharp-moment}]
The path supremum is invariant under the $J_1$ time changes, so Theorem \ref{thm:main} gives $\|U^n\|_\infty^p\Rightarrow\|U\|_\infty^p$.
To pass to expectations, take $p<q<\alpha$. Proposition \ref{prop:unif-mom} yields
\begin{equation*}
\sup_n\mE[\|U^n\|_\infty^p\mathbf1_{\{\|U^n\|_\infty>K\}}]
\leqslant K^{p-q}\sup_n\mE\|U^n\|_\infty^q\longrightarrow0,\qquad K\to\infty.
\end{equation*}
Uniform integrability and $U^n=u_n(X^n-X)$ now give
\begin{equation*}
u_n^p\mE\|X^n-X\|_\infty^p=\mE\|U^n\|_\infty^p
\longrightarrow \mE\|U\|_\infty^p.
\end{equation*}

For strict positivity, assume $\sigma_0(x_0)\sigma_0'(x_0)\ne0$. As $s\downarrow0$, $X_{s-}\to x_0$ almost surely and $\mu_{X_s}\to\delta_{x_0}$ in $\mathbb W_1$. Thus $\sigma_s(X_{s-})\sigma_s'(X_{s-})$ stays nonzero on a short interval after zero. Independence and the absence of fixed jumps of $V$ imply that $L$ and $V$ have no common jumps almost surely, by conditioning on the countable jump times of $L$. Also, for every rational $r\in(0,T]$,
\begin{equation*}
\mP\bigl(V\text{ has no jump in }(0,r]\bigr)
=\lim_{\varepsilon\downarrow0}
\mathrm e^{-r\nu_V(\{|z|>\varepsilon\})}=0.
\end{equation*}
Almost surely the conclusion holds simultaneously for all such rationals. Choose a rational $r>0$ in the interval above and a jump time $s\in(0,r]$ of $V$. Since $\Delta L_s=0$, \eqref{eq:reduced-limit} gives
\begin{equation*}
\Delta U_s=-\sigma_s(X_{s-})\sigma_s'(X_{s-})\Delta V_s\ne0 .
\end{equation*}
Thus $\mP(\|U\|_\infty>0)=1$, which gives strict positivity of the limit in \eqref{eq:path-moment-asymptotic}.
\end{proof}

\begin{proof}[Proof of Corollary \ref{cor:feedback-bias}]
Every subsequence has a further one along which $(L,M^n,N^n,U^n)$ converges by Proposition \ref{prop:error-tightness}. Theorem \ref{thm:main} identifies the law of $(X,Z)$ with that of $(X,U)$, since $X$ is a function of $L$. Corollary \ref{cor:meanfield-vanishes} therefore gives $F_Z^f=0$. Along this further subsequence, \eqref{eq:meanfield-Lr} gives the asserted coefficient convergence for every $1\leqslant r<\infty$. The drift estimate in Lemma \ref{lem:meanfield-passage} and \eqref{eq:meanfield-stochastic-passage} then give $\mE\|D^n\|_\infty^r\to0$ and $\mE\|S^n\|_\infty^p\to0$ for $1\leqslant p<\alpha$. The subsequence principle gives all three conclusions along the full sequence.

For the weak error at fixed $t$, take $\varphi\in\mathbf C_b^2(\mR)$ and $1<q<\alpha$. Taylor's formula and boundedness of $\varphi'$ give
\begin{equation*}
|\varphi(x+z)-\varphi(x)-\varphi'(x)z|
\leqslant C_\varphi(|z|^2\wedge|z|)\leqslant C_\varphi|z|^q.
\end{equation*}
For deterministic $t$, the Taylor bound and Proposition \ref{prop:unif-mom} give
\begin{equation*}
\left|u_n[\mE\varphi(X_t^n)-\mE\varphi(X_t)]-\mE[\varphi'(X_t)U_t^n]\right|
\leqslant C_\varphi u_n\mE|X_t^n-X_t|^q
=C_\varphi u_n^{1-q}\mE|U_t^n|^q\longrightarrow0 .
\end{equation*}
The processes $X,U$ have no fixed jumps by Lemma \ref{lem:wellposed}(2)
and by \eqref{eq:uniform-time-continuity-error} with Lemma \ref{lem:moving-evaluation}. Theorem \ref{thm:main} and Lemma \ref{lem:fixed-marginal} then give $(X_t,U_t^n)\Longrightarrow(X_t,U_t)$. Boundedness of $\varphi'$ and the $q$th moment bound give uniform integrability of $\varphi'(X_t)U_t^n$. Since $\varphi'(X_t)$ is a
function of $L$, conditional centering gives
\begin{equation*}
\mE[\varphi'(X_t)U_t^n]\longrightarrow\mE[\varphi'(X_t)U_t]
=\mE\left[\varphi'(X_t)\mE[U_t\mid\sigma(L_s:0\leqslant s\leqslant T)]\right]=0 .
\end{equation*}
Together with the Taylor bound this proves \eqref{eq:smooth-weak-bias}.
\end{proof}

\begin{proof}[Proof of Corollary \ref{cor:exact-law-comparison}]
Fix $1\leqslant p<\alpha$. The grid recursion in Lemma \ref{lem:wellposed} defines $Y^n$ using the exact laws $\mu_{X_{t_n}}$. The proof of part (1) gives the same path moment bound for $Y^n$, using boundedness of $\sigma$ and the growth estimate $\mE|b(Y_{s_n}^n,\mu_{X_{s_n}})|^p\leqslant C_p(1+\mE|Y_{s_n}^n|^p)$. Subtract the schemes driven by $L$ and split each coefficient difference at $(X_{s_n}^n,\mu_{X_{s_n}})$ to get
\begin{equation*}
\begin{aligned}
u_n(X_t^n-Y_t^n)
={}&D_t^n+S_t^n
+u_n\int_0^t[b(X_{s_n}^n,\mu_{X_{s_n}})-b(Y_{s_n}^n,\mu_{X_{s_n}})]\dif s\\
&+u_n\int_0^t[\sigma(X_{s_n}^n,\mu_{X_{s_n}})-\sigma(Y_{s_n}^n,\mu_{X_{s_n}})]\dif L_s,
\end{aligned}
\end{equation*}
where $D^n$ and $S^n$ are defined in Corollary \ref{cor:feedback-bias}. On each interval $[S,R]$ with grid endpoints, $s_n\in[S,R]$ for $S<s\leqslant R$. Each coefficient difference is at most $C|X_{s_n}^n-Y_{s_n}^n|$ in absolute value. The increments of $D^n,S^n$ are bounded by twice their path suprema. Starting at $S$, apply H\"older's inequality to the drift integral and Lemma \ref{lem:stable-maximal} to the jump integral to obtain
\begin{equation*}
\begin{aligned}
u_n^p\mE\sup_{S\leqslant t\leqslant R}|X_t^n-Y_t^n|^p
&\leqslant C_p\bigl(u_n^p\mE|X_S^n-Y_S^n|^p
+\mE\|D^n\|_\infty^p+\mE\|S^n\|_\infty^p\bigr)\\
&\quad+C_p\bigl((R-S)^p+(R-S)^{p/\alpha}\bigr)
u_n^p\mE\sup_{S\leqslant t\leqslant R}|X_t^n-Y_t^n|^p .
\end{aligned}
\end{equation*}
Choose $h_0>0$, independent of $n$, so the last coefficient is at most $1/2$ on intervals of length at most $h_0$. For large $n$, use the partition $0=r_0<\cdots<r_m=T$ in Lemma \ref{lem:wellposed}(4), with $m\leqslant 1+2T/h_0$. Set
\begin{equation*}
B_j^n=u_n^p\mE\sup_{r_j\leqslant t\leqslant r_{j+1}}|X_t^n-Y_t^n|^p,
\qquad
\varepsilon_n=\mE\|D^n\|_\infty^p+\mE\|S^n\|_\infty^p.
\end{equation*}
Absorbing the last term and using $X_0^n=Y_0^n$ and
$u_n^p\mE|X_{r_j}^n-Y_{r_j}^n|^p\leqslant B_{j-1}^n$ for $j\geqslant1$ gives
\begin{equation*}
B_0^n\leqslant C_p\varepsilon_n,\qquad
B_j^n\leqslant C_p(B_{j-1}^n+\varepsilon_n),\quad 1\leqslant j<m.
\end{equation*}
Since $m$ is uniformly bounded in $n$, induction gives
$B_j^n\leqslant C_{p,T}\varepsilon_n$, and hence
\begin{equation*}
u_n^p\mE\|X^n-Y^n\|_\infty^p
\leqslant\sum_{j=0}^{m-1}B_j^n
\leqslant C_{p,T}\varepsilon_n
\longrightarrow0.
\end{equation*}
The limit follows from Corollary \ref{cor:feedback-bias} with $r=p$ for the drift. Since $u_n(Y^n-X)=U^n-u_n(X^n-Y^n)$ and $u_n\|X^n-Y^n\|_\infty\to0$ in probability, the joint limit with $L$ transfers to $u_n(Y^n-X)$. Apply Lemma \ref{lem:fixed-marginal} for stable convergence. For the moment limit, Minkowski's inequality and Corollary \ref{cor:sharp-moment} give
\begin{equation*}
\begin{aligned}
&\left|u_n(\mE\|Y^n-X\|_\infty^p)^{1/p}-(\mE\|U\|_\infty^p)^{1/p}\right|\\
&\quad\leqslant u_n(\mE\|X^n-Y^n\|_\infty^p)^{1/p}
+\left|(\mE\|U^n\|_\infty^p)^{1/p}-(\mE\|U\|_\infty^p)^{1/p}\right|
\longrightarrow0 .
\end{aligned}
\end{equation*}
Taking $p$th powers gives the claimed moment asymptotic for $Y^n-X$.
\end{proof}

\noindent{\bf Acknowledgements}.
T. Wang is supported by the National Natural Science Foundation of China (Grant No.\ 12401178).

\raggedbottom
\makeatletter
\renewcommand{\@bibtitlestyle}{%
  \@startsection{section}{1}{\z@}{6pt}{4pt}%
    {\normalfont\scshape\centering}*{\refname}}
\makeatother
\bibliographystyle{plainurl}
\bibliography{EM_ERROR_LEVY}

\end{document}